\documentclass[11pt]{article}
\usepackage{bbm}
\usepackage{amssymb}
\usepackage{mathrsfs}
\usepackage{tikz}
\usepackage{mathrsfs}
\usepackage{amssymb}
\usepackage{amsmath}
\usepackage{bm}
\usepackage{hyperref}
\usepackage{graphicx}
\usepackage{pst-poly}  
\usepackage{pst-plot}  
\usepackage{rotating}

\newtheorem{thm}{Theorem}[section]
\newtheorem{lem}[thm]{Lemma}
\newtheorem{Def}[thm]{Definition}
\newtheorem{cor}[thm]{Corollary}

\newtheorem{con}[thm]{Conjecture}
\newtheorem{exa}[thm]{Example}

\newenvironment{pf}[1][Proof]{\noindent\textbf{#1.} }{\hfill\rule{2mm}{2mm}}

\makeatletter \@addtoreset{equation}{section} \makeatother

\begin{document}

\title
{Induced subgraphs of product graphs and a generalization of Huang's theorem}

\author
{Zhen-Mu Hong\\
{\footnotesize School of Finance, Anhui University of Finance \& Economics, Bengbu, Anhui 233030, China} \\
\\
Hong-Jian Lai\thanks{Corresponding author. E-mail addresses: zmhong@mail.ustc.edu.cn
(Z.-M. Hong), hjlai@math.wvu.edu (H.-J. Lai), jl0068@mix.wvu.edu (J.-B. Liu).}, Jian-Bing Liu\\
{\footnotesize Department of Mathematics, West Virginia University, Morgantown, WV 26506, USA}}

\date{}
\maketitle


\begin{center} {\bf Abstract} \end{center}

Recently, Huang showed that every $(2^{n-1}+1)$-vertex induced subgraph of the $n$-dimensional hypercube
has maximum degree at least $\sqrt{n}$ in [Annals of Mathematics, 190 (2019), 949--955].
In this paper, we discuss the induced subgraphs of Cartesian product graphs and semi-strong product graphs to generalize Huang's result.
Let $\Gamma_1$ be a connected signed bipartite graph of order $n$ and $\Gamma_2$
be a connected signed graph of order $m$. By defining two kinds of signed product
of $\Gamma_1$ and $\Gamma_2$, denoted by $\Gamma_1\widetilde{\Box}\Gamma_2$ and
$\Gamma_1\widetilde{\bowtie} \Gamma_2$, we show that if $\Gamma_1$ and $\Gamma_2$
have exactly two distinct adjacency eigenvalues $\pm\theta_1$ and $\pm\theta_2$ respectively,
then every $(\frac{1}{2}mn+1)$-vertex induced subgraph of
$\Gamma_1\widetilde{\Box}\Gamma_2$ (resp. $\Gamma_1\widetilde{\bowtie} \Gamma_2$)
has maximum degree at least $\sqrt{\theta_1^2+\theta_2^2}$ (resp. $\sqrt{(\theta_1^2+1)\theta_2^2}$).
Moreover, we discuss the eigenvalues of $\Gamma_1\widetilde{\Box} \Gamma_2$ and $\Gamma_1\widetilde{\bowtie} \Gamma_2$
and obtain a sufficient and necessary condition such that the spectrum of
$\Gamma_1\widetilde{\Box}\Gamma_2$ and $\Gamma_1\widetilde{\bowtie}\Gamma_2$
are symmetric, from which we obtain more general results on maximum degree of the induced subgraphs.

\vskip10pt




\section{Introduction}

Let $Q_n$ be the $n$-dimensional hypercube, whose vertex set consists of vectors
in $\{0, 1\}^n$, and two vectors are adjacent if they differ in exactly one coordinate.
For a simple and undirected graph $G=(V,E)$, we use $\Delta(G)$ to denote the maximum degree of $G$.
The adjacency matrix of $G$ is defined to be a $(0,1)$-matrix
$A(G)=(a_{ij})$, where $a_{ij}=1$ if $v_i$ and $v_j$ are adjacent, and
$a_{ij}=0$ otherwise.

Recently, Huang \cite{Huan19} constructed a signed adjacency matrix of $Q_n$
with exactly two distinct eigenvalues $\pm \sqrt{n}$. Using eigenvalue interlacing, Huang
proceeded to prove that the spectral radius (and so, the maximum degree) of any $(2^{n-1} + 1)$-vertex induced
subgraph of $Q_n$, is at least $\sqrt{n}$. Combing this with the combinatorial equivalent
formulation discovered by Gotsman and Linial \cite{GoLi92},
Huang confirmed the Sensitivity Conjecture \cite{NiSz92} from theoretical computer science.
The main contribution of Huang is the following theorem.

\begin{thm}{(Huang \cite{Huan19})}\label{thm-Huang}
For every integer $n\geq 1$, let $H$ be an arbitrary $(2^{n-1}+1)$-vertex induced
subgraph of $Q_n$, then $\Delta(H)\geq \sqrt{n}$.
\end{thm}

The bound $\sqrt{n}$ (or more precisely, $\lceil\sqrt{n}\rceil$)
is sharp, as shown by Chung, Furedi, Graham, and Seymour \cite{CFGS88} in 1988.
Tao \cite{Tao19} also gave a great expository of Huang's work on his blog after Huang
announced the proof of the Sensitivity Conjecture.

Denote the Cartesian product of two graphs $G$ and $H$ by $G \Box H$.
It is known that the hypercube $Q_n$ can be constructed iteratively by Cartesian product, that is,
$Q_1=K_2$ and for $n\geq 2$, $Q_n=Q_1\Box Q_{n-1}$. Motivated by this fact,
in this paper, we generalize Huang's theorem to
Cartesian product graphs and semi-strong product graphs.
We introduce some necessary notations in the following.

A {\it signed graph} $\Gamma = (G, \sigma)$ is a graph $G = (V, E)$,
together with a sign function $\sigma : E \rightarrow \{+1,-1\}$ assigning
a positive or negative sign to each edge. An edge $e$ is positive if $\sigma(e) = 1$
and negative if $\sigma(e) = -1$. The unsigned graph $G$ is said to
be the {\em underlying graph} of $\Gamma$, while $\sigma$ is called the
{\em signature} of $G$. If each edge of $\Gamma$ is positive (resp. negative),
then $\Gamma$ is denoted by $\Gamma= (G,+)$ (resp. $\Gamma= (G,-)$).
A signed graph is connected if its underlying graph is connected.

The adjacency matrix of $\Gamma=(G,\sigma)$ is denoted by $A(\Gamma)=(a^{\sigma}_{ij})$, where
$a^{\sigma}_{ij} = \sigma(v_iv_j)$, if $v_i$ and $v_j$ are adjacent, and $a^{\sigma}_{ij}=0$ otherwise.
As $G$ is simple and undirected, the adjacency matrix $A(\Gamma)$ is a symmetric $(-1, 0, +1)$-matrix,
and $A(\Gamma)=A(G)$ if $\Gamma= (G,+)$, $A(\Gamma)=-A(G)$ if $\Gamma= (G,-)$.
Let $\lambda_1(\Gamma)\geq \lambda_2(\Gamma)\geq\cdots\geq\lambda_n(\Gamma)$ denote the eigenvalues of
$A(\Gamma)$, which are all real since $A(\Gamma)$ is real and symmetric.
If $\Gamma$ contains at least one edge, then $\lambda_1(\Gamma) > 0 > \lambda_n(\Gamma)$
since the trace of $A(\Gamma)$ is $0$. In general, the largest eigenvalue $\lambda_1(\Gamma)$ may not be equal to the spectral radius
$\rho(\Gamma) = \max\{|\lambda_i(\Gamma)| : 1 \leq i \leq n\} = \max\{\lambda_1(\Gamma),-\lambda_n(\Gamma)\}$ because
the Perron-Frobenius Theorem is valid only for nonnegative matrices.
The eigenvalues of the adjacency matrix of signed graph $\Gamma$ are called adjacency eigenvalues of $\Gamma$.
The spectrum of $A(\Gamma)$ is called the (adjacency) spectrum of $\Gamma$
and $A(\Gamma)$ is also called a signed adjacency matrix of $G$.
The spectrum of $\Gamma$ is {\it symmetric} if its adjacency eigenvalues
are symmetric with respect to the origin.
In this paper, all eigenvalues considered are adjacency eigenvalues.

For basic results in the theory of signed graphs, the reader is referred to Zaslavsky \cite{Zasl82}.
Recently, the spectra of signed graphs have attracted much attention,
as found in \cite{ABDN18,BCKW19,BiLi06,Gall16,GhFa17,GeHZ11,HoTW19,Rame15,Rame19,Stan20,Zasl10}, among others.
In \cite{BCKW19}, the authors surveyed some general results on the adjacency
spectra of signed graphs and proposed some spectral problems which are inspired by the
spectral theory of unsigned graphs. In particular, the signed graphs with exactly
two distinct eigenvalues have been greatly investigated in recent years,
see \cite{GhFa17,HoTW19,McSm07,Rame15,Rame19,Stan20}.
In \cite{HoTW19}, Hou et al. characterized all simple connected
signed graphs with maximum degree at most 4 and with just two distinct adjacency eigenvalues.
In this paper, we construct signed graphs with exactly two distinct eigenvalues
by two kinds of graph products, which generalizes Huang's result on the induced subgraph of the hypercube.

The {\it Kronecker product} $A\otimes B$ of matrices $A = (a_{ij})_{m\times n}$
and $B = (b_{ij})_{p\times q}$ is the $mp \times nq$ matrix obtained from $A$ by replacing each
element $a_{ij}$ with the block $a_{ij}B$.
Therefore the entries of $A\otimes B$ consist of all the $mnpq$ possible products of an entry
of $A$ with an entry of $B$. For matrices $A, B, C$ and $D$, we have
$(A \otimes B) \cdot (C \otimes D) = AC \otimes BD$ whenever the products $AC$ and $BD$ exist.
Note that, $(A \otimes B)^T=A^T\otimes B^T$.

The {\it Cartesian product} of two graphs $G_1$ and $G_2$ is a graph, denoted by $G_1 \Box G_2$, whose vertex set is
$V(G_1)\times V(G_2)$ and two vertices $(u_1, u_2)$ and $(v_1, v_2)$ being adjacent in $G_1 \Box G_2$
if and only if either $u_1=v_1$ and $u_2v_2\in E(G_2)$, or $u_1v_1 \in E(G_1)$ and $u_2=v_2$.
The {\it direct product} (or {\it Kronecker product}) of two graphs $G_1$ and $G_2$ is a graph, denoted by $G_1 \times G_2$,
whose vertex set is $V(G_1)\times V(G_2)$, and two vertices
$(u_1, u_2)$ and $(v_1, v_2)$ being adjacent to each other in $G_1 \times G_2$
if and only if both $u_1v_1 \in E(G_1)$ and $u_2v_2\in E(G_2)$.
The {\it semi-strong product} (or {\it strong tensor product} \cite{GaRW76})
of two graphs $G_1$ and $G_2$ is a graph, denoted by $G_1 \bowtie G_2$,
whose vertex set is $V(G_1)\times V(G_2)$, and two vertices
$(u_1, u_2)$ and $(v_1, v_2)$ being adjacent to each other in $G_1 \bowtie G_2$
if and only if either $u_1v_1 \in E(G_1)$ and $u_2v_2\in E(G_2)$, or $u_1=v_1$ and $u_2v_2\in E(G_2)$.
Then, by the definitions, the adjacency matrices of
$G_1\Box G_2$, $G_1\times G_2$ and $G_1\bowtie G_2$ are
$A(G_1\Box G_2) =A(G_1)\otimes I_m+I_n\otimes A(G_2)$,
$A(G_1\times G_2) =A(G_1)\otimes A(G_2)$ and $A(G_1\bowtie G_2) =A(G_1)\otimes A(G_2)+I_n\otimes A(G_2)$, respectively,
where $n=|V(G_1)|$, $m=|V(G_2)|$ and $I_n$ is the identity matrix of order $n$.

Let $\Gamma_1=(G_1,\sigma_1)$ be a connected signed bipartite graph of order $n$ with bipartition $(V_1,V_2)$,
where $|V_1|=s$ and $|V_2|=n-s$, and $\Gamma_2=(G_2,\sigma_2)$ be a connected signed graph of order $m$.
With suitable labeling of vertices, the adjacency matrix of $\Gamma_1$ can be represented as
$$
A(\Gamma_1)=\left[
\begin{array}{cc}
O_{s} & P \\
P^T & O_{n-s}
\end{array}
\right].
$$

The {\it signed Cartesian product} of signed bipartite graph $\Gamma_1$ and signed graph $\Gamma_2$,
denoted by $\Gamma_1\widetilde{\Box} \Gamma_2$, is the signed graph with adjacency matrix
\begin{equation}\label{e1.1}
A(\Gamma_1\widetilde{\Box} \Gamma_2) =  A(\Gamma_1) \otimes I_m +
\left[
\begin{array}{cc}
I_s & O \\
O & -I_{n-s}
\end{array}
\right] \otimes A(\Gamma_2)=\left[
\begin{matrix}
I_s\otimes A(\Gamma_2) & P\otimes I_m \\
P^T\otimes I_m & -I_{n-s}\otimes A(\Gamma_2)
\end{matrix}
\right].
\end{equation}

The {\it signed semi-strong product} of signed bipartite graph $\Gamma_1$ and signed graph $\Gamma_2$,
denoted by $\Gamma_1 \widetilde{\bowtie} \Gamma_2$, is the signed graph with adjacency matrix
\begin{equation}\label{e1.2}
A(\Gamma_1\widetilde{\bowtie} \Gamma_2) = A(\Gamma_1)\otimes A(\Gamma_2) +
\left[
\begin{array}{cc}
I_s & O \\
O & -I_{n-s}
\end{array}
\right] \otimes A(\Gamma_2)=\left[
\begin{array}{cc}
I_s & P \\
P^T & -I_{n-s}
\end{array}
\right] \otimes A(\Gamma_2).
\end{equation}

As a generalization of Theorem \ref{thm-Huang}, we have the following theorem.

\begin{thm}\label{thm-DeltaH-2}
Let $\Gamma_1=(G_1,\sigma_1)$ be a signed bipartite graph of order $n$
with exactly two distinct eigenvalues $\pm \theta_1$ and
$\Gamma_2=(G_2,\sigma_2)$ be a signed graph of order $m$
with exactly two distinct eigenvalues $\pm \theta_2$.
If $H$ and $H'$ are arbitrary $(\frac{mn}{2}+1)$-vertex induced
subgraphs of $\Gamma_1\widetilde{\Box} \Gamma_2$ and
$\Gamma_1\widetilde{\bowtie} \Gamma_2$ respectively, then
$$
\Delta(H)\geq \sqrt{\theta_1^2+\theta_2^2}, \ \
\Delta(H')\geq \sqrt{(\theta_1^2+1)\theta_2^2}.
$$
\end{thm}

A direct proof of Theorem \ref{thm-DeltaH-2} is presented in Section 2.
In fact, from the proof we see that $\Gamma_1$ and $\Gamma_2$ in Theorem \ref{thm-DeltaH-2} are regular.
For more general graphs, we can obtain the following theorem.
A (signed) bipartite graph with bipartition $(V_1,V_2)$ is called {\it balanced}
if $|V_1|=|V_2|$.

\begin{thm}\label{thm-Delta-2-min}
Let $\Gamma_1=(G_1,\sigma_1)$ be a signed bipartite graph of order $n$
and $\Gamma_2=(G_2,\sigma_2)$ be a signed graph of order $m$,
and let $\lambda^2$ and $\mu^2$ be the minimum eigenvalues of $A(\Gamma_1)^2$ and $A(\Gamma_2)^2$, respectively.
Let $H$ and $H'$ be any $(\lfloor\frac{mn}{2}\rfloor+1)$-vertex induced
subgraph of $\Gamma_1\widetilde{\Box} \Gamma_2$ and
$\Gamma_1\widetilde{\bowtie} \Gamma_2$, respectively.
If $\Gamma_1$ is a balanced bipartite graph or the spectrum of $\Gamma_2$ is symmetric, then
$$
\Delta(H)\geq \sqrt{\lambda^2+\mu^2}, \ \  \Delta(H')\geq \sqrt{(\lambda^2+1)\mu^2}.
$$
\end{thm}

In Section 3, we display some preliminaries and examples.
In Section 4, we give a characterization of the eigenvalues of
$\Gamma_1\widetilde{\Box} \Gamma_2$ and $\Gamma_1\widetilde{\bowtie} \Gamma_2$
and obtain a sufficient and necessary condition such that the spectrum of
$\Gamma_1\widetilde{\Box}\Gamma_2$ and $\Gamma_1\widetilde{\bowtie}\Gamma_2$
are symmetric. In Section 5, we present the proof of Theorem \ref{thm-Delta-2-min}
and generalize the signed Cartesian product and signed semi-strong product of
two signed graphs to the products of $n$ signed graphs.
In the last section, we give some concluding remarks.

\section{A direct proof of Theorem \ref{thm-DeltaH-2}}

Using the idea that Shalev Ben-David contributed on July 3, 2019 to Scott Aaronson's blog,
Knuth \cite{knot19} gave a direct and nice proof of Huang's theorem in one page. Here,
arising from their ideas, we give a direct proof of Theorem \ref{thm-DeltaH-2}.

\bigskip
\noindent{\bf Proof of Theorem \ref{thm-DeltaH-2}.}
For simplicity, let $A_1:=A(\Gamma_1)$ and $A_2:=A(\Gamma_2)$.
Since $\Gamma_i$ has exactly two distinct eigenvalues $\pm\theta_i~(\neq 0)$ for $i=1,2$, we have
each eigenvalue of $A_i^2$ equals to $\theta_i^2$ and so there exist
orthogonal matrices $Q_1$ and $Q_2$ such that $A_1^2=Q_1(\theta_1^2I_n)Q_1^T=\theta_1^2 I_n$
and $A_2^2=Q_2(\theta_2^2I_m)Q_2^T=\theta_2^2I_m$.
The diagonal entries of $A_i^2$ are the degrees of vertices in $\Gamma_i$, so
$\Gamma_i$ is a $\theta_i^2$-regular graph for $i=1,2$.
Moreover, $|V_1|=s=\frac{n}{2}$ and $PP^T=P^TP=\theta_1^2I_{n/2}$.

(a) Let $\mathcal {A}:=A(\Gamma_1\widetilde{\Box} \Gamma_2)$ and define
$$
\mathcal{B}=\left[
\begin{matrix}
P\otimes (A_2+\sqrt{\theta_1^2+\theta_2^2}I_m) \\
\theta_1^2I_{n/2}\otimes I_m
\end{matrix}
\right]
$$
to be an $mn\times \frac{mn}{2}$ matrix.
Since $\theta_1\neq 0$, the rank of $\mathcal{B}$ is $\frac{mn}{2}$, and we have
\begin{align}
\mathcal{A}\cdot \mathcal{B} 
= &\left[
\begin{matrix}
I_{n/2}\otimes A_2 & P\otimes I_m \\
P^T\otimes I_m & -I_{n/2}\otimes A_2
\end{matrix}
\right] \cdot \left[
\begin{matrix}
P\otimes (A_2+\sqrt{\theta_1^2+\theta_2^2}I_m) \\
\theta_1^2I_{n/2}\otimes I_m
\end{matrix}
\right] \nonumber\\ 
= & \left[
\begin{matrix}
P\otimes (A_2^2+\sqrt{\theta_1^2+\theta_2^2}A_2+\theta_1^2I_m) \\
P^TP\otimes (A_2+\sqrt{\theta_1^2+\theta_2^2}I_m)-\theta_1^2I_{n/2}\otimes A_2
\end{matrix}
\right] \nonumber\\ 
= & \left[
\begin{matrix}
P\otimes (\theta_2^2I_m+\sqrt{\theta_1^2+\theta_2^2}A_2+\theta_1^2I_m) \\
\theta_1^2I_{n/2}\otimes (A_2+\sqrt{\theta_1^2+\theta_2^2}I_m)-\theta_1^2I_{n/2}\otimes A_2
\end{matrix}
\right] \nonumber\\ 
= & \sqrt{\theta_1^2+\theta_2^2} \left[
\begin{matrix}
P\otimes (A_2+\sqrt{\theta_1^2+\theta_2^2}I_m) \\
\theta_1^2I_{n/2}\otimes I_m
\end{matrix}
\right]
=  \sqrt{\theta_1^2+\theta_2^2} \mathcal{B}. \nonumber
\end{align}

Let $H$ be an arbitrary $(\frac{mn}{2}+1)$-vertex induced subgraph of $\Gamma_1\widetilde{\Box} \Gamma_2$.
Suppose $\mathcal{B}^*$ is the $(\frac{mn}{2}-1)\times \frac{mn}{2}$
submatrix of $\mathcal{B}$ whose rows corresponding to vertices not in $H$.
Then there exists a unit $\frac{mn}{2}\times 1$ vectors $x$ such that $\mathcal{B}^*x=0$,
since $\mathcal{B}^*x=0$ is a homogeneous system of $\frac{mn}{2}-1$ linear equations with $\frac{mn}{2}$ variables.
As ${\rm rank}(\mathcal{B})=\frac{mn}{2}$, $y=\mathcal{B}x$ is an $mn\times 1$
nonzero vector such that $y_v=0$ for any vertex $v\not\in H$, and
$\mathcal{A}y =\sqrt{\theta_1^2+\theta_2^2} y$.

Let $u$ be a vertex such that $|y_u|=\max\{|y_1|,\dots,|y_{mn}|\}$. Then $|y_u|>0$, $u\in V(H)$ and
$$
\sqrt{\theta_1^2+\theta_2^2} |y_u|=|(\mathcal{A}y)_{u}|=\left|\sum_{v=1}^{mn}\mathcal{A}_{uv}y_v\right|
=\left|\sum_{v\in H}\mathcal{A}_{uv}y_v\right|\leq \sum_{v\in H}|\mathcal{A}_{uv}||y_u|
\leq \Delta(H)|y_u|.
$$
Therefore, $\Delta(H)\geq \sqrt{\theta_1^2+\theta_2^2}$.

(b) Let $\mathcal{A}:=A(\Gamma_1\widetilde{\bowtie} \Gamma_2)$
and define
$$
\mathcal{B}=\left[
\begin{matrix}
P\otimes(\sqrt{\theta_1^2+1}A_2+\theta_2I_m) \\
\theta_1^2\theta_2I_{n/2}\otimes I_m
\end{matrix}
\right]
$$
to be an $mn\times \frac{mn}{2}$ matrix.
Since $\theta_1\neq 0$ and $\theta_2\neq 0$,
the rank of $\mathcal{B}$ is $\frac{mn}{2}$, and we have
\begin{align}
\mathcal{A}\cdot \mathcal{B} 
= & \left[
\begin{matrix}
I_{n/2}\otimes A_2 & P\otimes A_2  \\
P^T\otimes A_2 & -I_{n/2}\otimes A_2
\end{matrix}
\right] \cdot \left[
\begin{matrix}
P\otimes(\sqrt{\theta_1^2+1}A_2+\theta_2I_m) \\
\theta_1^2\theta_2I_{n/2}\otimes I_m
\end{matrix}
\right] \nonumber\\ 
= & \left[
\begin{matrix}
P\otimes(\sqrt{\theta_1^2+1}A_2^2+\theta_2A_2+\theta_1^2\theta_2A_2) \\
P^TP\otimes (\sqrt{\theta_1^2+1}A_2^2+\theta_2A_2)-\theta_1^2\theta_2I_{n/2}\otimes A_2
\end{matrix}
\right] \nonumber\\ 
= & \left[
\begin{matrix}
P\otimes(\sqrt{\theta_1^2+1}\theta_2^2I_m+\theta_2(\theta_1^2+1)A_2) \\
\theta_1^2I_{n/2}\otimes (\sqrt{\theta_1^2+1}\theta_2^2I_m+\theta_2A_2)-\theta_1^2\theta_2I_{n/2}\otimes A_2
\end{matrix}
\right] \nonumber\\ 
= & \sqrt{(\theta_1^2+1)\theta_2^2} \left[
\begin{matrix}
P\otimes(\sqrt{\theta_1^2+1}A_2+\theta_2I_m) \\
\theta_1^2\theta_2I_{n/2}\otimes I_m
\end{matrix}
\right]
= \sqrt{(\theta_1^2+1)\theta_2^2} \mathcal{B}. \nonumber
\end{align}

Let $H'$ be an arbitrary $(\frac{mn}{2}+1)$-vertex induced subgraph of $\Gamma_1\widetilde{\bowtie} \Gamma_2$.
Suppose $\mathcal{B}^*$ is the $(\frac{mn}{2}-1)\times \frac{mn}{2}$ submatrix of $\mathcal{B}$
whose rows corresponding to vertices not in $H'$.
Then there exists a unit $\frac{mn}{2}\times 1$ vector $x$ such that $\mathcal{B}^*x=0$,
since $\mathcal{B}^*x=0$ is a homogeneous system of $\frac{mn}{2}-1$ linear equations with $\frac{mn}{2}$ variables.
As ${\rm rank}(\mathcal{B})=\frac{mn}{2}$, $y=\mathcal{B}x$ is an $mn\times 1$ nonzero vector
such that $y_v=0$ for any vertex $v\not\in H'$, and
$\mathcal{A}y =\sqrt{(\theta_1^2+1)\theta_2^2} y$.

Let $u$ be a vertex such that $|y_u|=\max\{|y_1|,\dots,|y_{mn}|\}$. Then $|y_u|>0$, $u\in V(H')$ and
$$
\sqrt{(\theta_1^2+1)\theta_2^2} |y_u|=|(\mathcal{A}y)_{u}|=\left|\sum_{v=1}^{mn}\mathcal{A}_{uv}y_v\right|
=\left|\sum_{v\in H'}\mathcal{A}_{uv}y_v\right|\leq \sum_{v\in H'}|\mathcal{A}_{uv}||y_u|
\leq \Delta(H')|y_u|.
$$
Therefore, $\Delta(H')\geq \sqrt{(\theta_1^2+1)\theta_2^2}$.  \hfill\rule{2mm}{2mm}

\section{Preliminaries}

In this section, we present some useful lemmas and examples.

\begin{lem}{(Hammack et al. \cite{HaIK11})}\label{lem-HaIK11}
Let $G_1$ and $G_2$ be nontrivial graphs. Then

(i) $G_1\Box G_2$ is connected if and only if $G_1$ and $G_2$ are connected,
and $G_1\Box G_2$ is bipartite if and only if $G_1$ and $G_2$ are bipartite;

(ii) $G_1\times G_2$ is connected if and only if $G_1$ and $G_2$ are connected and at most one of them is bipartite,
and $G_1\times G_2$ is bipartite if and only if at least one of $G_1$ and $G_2$ is bipartite.
\end{lem}

\begin{lem}{(Garman et al. \cite{GaRW76})}\label{lem-GaRW76}
Let $G_1$ and $G_2$ be nontrivial graphs. Then

(i) $G_1\bowtie G_2$ is connected if and only if $G_1$ and $G_2$ are connected;

(ii) $G_1\bowtie G_2$ is bipartite if and only if $G_2$ is bipartite;

(iii) The semi-strong product operation is neither associative nor commutative;

(iv) If $G_1$ is bipartite, then $G_1\bowtie K_2\cong G_1\Box K_2$.
\end{lem}

By Lemma \ref{lem-GaRW76} (iv), the following corollary can be obtained easily.

\begin{cor}{(Garman et al. \cite{GaRW76})}\label{cor-GaRW76}
(i) Let $G_1=K_2$, and for $n\geq 2$, $G_n=G_{n-1}\bowtie K_2$, then $G_n\cong Q_n$.

(ii) Let $G'_1=K_2$, and for $n\geq 2$, $G'_n=K_2\bowtie G'_{n-1}$, then $G'_n\cong K_{2^{n-1},2^{n-1}}$.
\end{cor}

\begin{pf}
By Lemma \ref{lem-GaRW76} (iv), $Q_{n-1}\bowtie K_2\cong Q_{n-1}\Box K_2=Q_n$. By induction, $G_n\cong Q_n$.
Let $V(K_2)=\{u,v\}$ and $(V_1,V_2)$ be the bipartition of $K_{2^{n-2},2^{n-2}}$.
Then there is an edge connecting any two vertices between $\{u,v\}\times V_1$ and $\{u,v\}\times V_2$
in $K_2\bowtie K_{2^{n-2},2^{n-2}}$. Hence, $K_2\bowtie K_{2^{n-2},2^{n-2}} = K_{2^{n-1},2^{n-1}}$.
By induction, $G'_n\cong K_{2^{n-1},2^{n-1}}$.
\end{pf}

\medskip
By the definitions of Cartesian product, direct product and semi-strong product of graphs,
we can define the product of signed graphs
$\Gamma_1$ and $\Gamma_2$ by their adjacency matrices. That is,
$A(\Gamma_1\Box \Gamma_2)=A(\Gamma_1)\otimes I_m+I_n\otimes A(\Gamma_2)$,
where $n=|V(\Gamma_1)|$ and $m=|V(\Gamma_2)|$,
$A(\Gamma_1\times \Gamma_2)=A(\Gamma_1)\otimes A(\Gamma_2)$, and
$A(\Gamma_1\bowtie \Gamma_2)=(A(\Gamma_1)+I_n)\otimes A(\Gamma_2)$.
If $X$ and $Y$ are eigenvectors of $A_1=A(\Gamma_1)$ and $A_2=A(\Gamma_2)$
corresponding to eigenvalues $\lambda$ and $\mu$, respectively, then direct computation yields the following.
\begin{align*}
&A(\Gamma_1\Box \Gamma_2)(X\otimes Y)     =(A_1\otimes I_m+I_n\otimes A_2)(X\otimes Y)=(\lambda+\mu) X\otimes Y,\\
&A(\Gamma_1\times \Gamma_2)(X\otimes Y)   =(A_1\otimes A_2)(X\otimes Y)=A_1X\otimes A_2Y=\lambda\mu X\otimes Y,\\
&A(\Gamma_1\bowtie \Gamma_2)(X\otimes Y)   =[(A_1+I_n)\otimes A_2](X\otimes Y)=(A_1+I_n)X\otimes A_2Y=(\lambda+1)\mu X\otimes Y.
\end{align*}

Thus, we can obtain the following theorem.

\begin{thm}\label{thm-Kronecker-eig}
If $\lambda_1\geq\lambda_2\geq\cdots\geq\lambda_n$ and $\mu_1\geq\mu_2\geq\cdots\geq\mu_m$
are the adjacency eigenvalues of the signed graphs $\Gamma_1$ and $\Gamma_2$, respectively,
then, for $i=1,2,\dots,n$ and $j=1,2,\dots,m$,

(i) (Germina et al. \cite{GeHZ11}) $\lambda_i+\mu_j$ are the adjacency eigenvalues of $\Gamma_1\Box \Gamma_2$;

(ii) (Germina et al. \cite{GeHZ11}) $\lambda_i\mu_j$ are the adjacency eigenvalues of $\Gamma_1\times \Gamma_2$;

(iii) $(\lambda_i+1)\mu_j$ are the adjacency eigenvalues of $\Gamma_1\bowtie \Gamma_2$.

\end{thm}

By Lemma \ref{lem-HaIK11} (ii) and Theorem \ref{thm-Kronecker-eig} (ii), we have the following result immediately.

\begin{cor}\label{cor-kron-eig}
For $i=1,2$, let $\Gamma_i=(G_i,\sigma_i)$ be a connected signed graph
with exactly two distinct eigenvalues $\pm\theta_i$, respectively.
If at least one of $G_1$ and $G_2$ is non-bipartite, then $\Gamma_1\times \Gamma_2$ is a connected
signed graph with exactly two distinct eigenvalues $\pm\theta_1\theta_2$.
\end{cor}

In the following, we introduce some known results and examples which can be used to construct
signed graphs with exactly two distinct eigenvalues. First we give some definitions.
A {\it weighing matrix} of order $n$ and weight $k$ is an $n\times n$ matrix $W=W(n,k)$ with
entries $0$, $+1$ and $-1$ such that $WW^T=W^T W=kI_{n}$.
A weighing matrix $W(n,n)$ is a Hadamard matrix $H_n$ of order $n$.
A {\it conference matrix} $C$ of order $n$ is an $n\times n$ matrix with $0$'s on
the diagonal, $+1$ or $-1$ in all other positions and with the property
$CC^T=(n-1)I_n$. Thus, a conference matrix of order $n$ is a weighing matrix of order $n$ and weight $n-1$,
and a permutation matrix of order $n$ is a weighing matrix of order $n$ and weight $1$.

\begin{lem}\label{lem-K-nn}
For $n\geq 1$, let
$$
H_2=\left[\begin{matrix}
   1 & 1 \\
   1 & -1 \\
\end{matrix}
\right],
H_{2^{n+1}} = H_2\otimes H_{2^{n}},
A_n=\left[\begin{matrix}
   0 & 1 \\
   1 & 0
\end{matrix}
\right]\otimes H_{2^n}.
$$
Then $A_n$ is a signed adjacency matrix of $K_{2^n,2^n}$
and its eigenvalues are $\pm\sqrt{2^n}$, each with multiplicity $2^n$.
\end{lem}

\begin{pf}
Since $H_{2^n}$ is a symmetric matrix with entries $\pm1$, $A_n$ is a signed adjacency matrix of $K_{2^n,2^n}$.
Note that $H_{2^n}$ is a Hadamard matrix of order $2^n$ with eigenvalues $\pm\sqrt{2^n}$.
By the property of Kronecker product, the eigenvalues of $A_n$ are $\pm\sqrt{2^n}$, each with multiplicity $2^n$.
\end{pf}

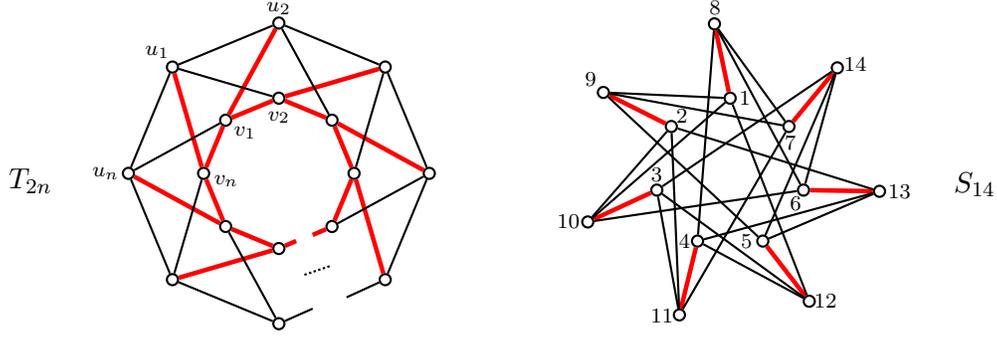
\begin{figure}[ht]
\begin{center}
\begin{pspicture}(-5.5,-2)(5.5,2)

\cnode(-4,0){2.5pt}{a1} \cnode(-3.7071,-0.7071){2.5pt}{h1}
\cnode(-3.7071,0.7071){2.5pt}{b1} \cnode(-3,-1){2.5pt}{g1}
\cnode(-3,1){2.5pt}{c1} \cnode(-2.2929,-0.7071){2.5pt}{f1}
\cnode(-2.2929,0.7071){2.5pt}{d1} \cnode(-2,0){2.5pt}{e1}

\rput(-6.3,-0.1){$T_{2n}$}
\rput(-3.7,-0.1){\scriptsize $v_n$}
\rput(-3.45,0.55){\scriptsize $v_1$}
\rput(-3,0.75){\scriptsize $v_2$}

\cnode(-5,0){2.5pt}{a} \cnode(-4.4142,-1.4142){2.5pt}{h}
\cnode(-4.4142,1.4142){2.5pt}{b} \cnode(-3,-2){2.5pt}{g}
\cnode(-3,2){2.5pt}{c} \cnode(-1.5858,-1.4142){2.5pt}{f}
\cnode(-1.5858,1.4142){2.5pt}{d} \cnode(-1,0){2.5pt}{e}

\rput(-5.3,0){\scriptsize $u_n$}
\rput(-4.6142,1.6142){\scriptsize $u_1$}
\rput(-3,2.2){\scriptsize $u_2$}

\ncline{a}{b}\ncline{b}{c}\ncline{c}{d}\ncline{d}{e}\ncline{e}{f}
\ncline{g}{h}\ncline{h}{a}
\ncline[linewidth=0.6mm,linecolor=red]{a1}{b1}\ncline[linewidth=0.6mm,linecolor=red]{b1}{c1}
\ncline[linewidth=0.6mm,linecolor=red]{c1}{d1}\ncline[linewidth=0.6mm,linecolor=red]{d1}{e1}
\ncline[linewidth=0.6mm,linecolor=red]{e1}{f1}
\ncline[linewidth=0.6mm,linecolor=red]{g1}{h1}\ncline[linewidth=0.6mm,linecolor=red]{h1}{a1}

\ncline[linewidth=0.6mm,linecolor=red]{a1}{b}\ncline[linewidth=0.6mm,linecolor=red]{b1}{c}
\ncline[linewidth=0.6mm,linecolor=red]{c1}{d}\ncline[linewidth=0.6mm,linecolor=red]{d1}{e}
\ncline[linewidth=0.6mm,linecolor=red]{e1}{f}
\ncline[linewidth=0.6mm,linecolor=red]{g1}{h}\ncline[linewidth=0.6mm,linecolor=red]{h1}{a}
\ncline{a}{b1}\ncline{b}{c1}\ncline{c}{d1}\ncline{d}{e1}
\ncline{e}{f1}
\ncline{g}{h1}\ncline{h}{a1}

\pnode(-2.77,-0.9){g1B}\pnode(-2.55,-0.82){f1B}\ncline[linewidth=0.6mm,linecolor=red]{g1}{g1B}\ncline[linewidth=0.6mm,linecolor=red]{f1}{f1B}
\pnode(-2.55,-1.8){gB}\pnode(-2.1,-1.65){fB}\ncline{g}{gB}\ncline{f}{fB}
\pnode(-2.66,-1.35){D1}\pnode(-2.325,-1.235){D2} \ncline[linestyle=dashed,dash=1pt 1pt]{D1}{D2}

\cnode(3,1){2.5pt}{u1}\cnode(2.2182,0.6235){2.5pt}{u2}\cnode(2.0251,-0.2225){2.5pt}{u3}
\cnode(2.5661,-0.9010){2.5pt}{u4}\cnode(3.4339,-0.9010){2.5pt}{u5}\cnode(3.9749,-0.2225){2.5pt}{u6}\cnode(3.7818,0.6235){2.5pt}{u7}

\rput(3.2,1){\scriptsize $1$}\rput(2.3582,0.7635){\scriptsize $2$}\rput(2.0251,-0.0225){\scriptsize $3$}
\rput(2.3961,-0.9010){\scriptsize $4$}\rput(3.2209,-0.9010){\scriptsize $5$}
\rput(3.8749,-0.4025){\scriptsize $6$}\rput(3.8018,0.4035){\scriptsize $7$}

\cnode(2.7909,1.9890){2.5pt}{v1}\cnode(1.3146,1.0767){2.5pt}{v2}\cnode(1.1073,-0.6464){2.5pt}{v3}
\cnode(2.3253,-1.8828){2.5pt}{v4}\cnode(4.0514,-1.7014){2.5pt}{v5}\cnode(4.9857,-0.2388){2.5pt}{v6}\cnode(4.4248,1.4036){2.5pt}{v7}

\rput(2.8,2.2){\scriptsize $8$}\rput(1.1463,1.2470){\scriptsize $9$}\rput(0.8501,-0.6450){\scriptsize $10$}
\rput(2.0922,-1.9019){\scriptsize $11$}\rput(4.2678,-1.7019){\scriptsize $12$}
\rput(5.2499,-0.2450){\scriptsize $13$}\rput(4.6637,1.4470){\scriptsize $14$}
\rput(6.2499,-0.1450){$S_{14}$}

\ncline[linewidth=0.6mm,linecolor=red]{u1}{v1}
\ncline[linewidth=0.6mm,linecolor=red]{u2}{v2}
\ncline[linewidth=0.6mm,linecolor=red]{u3}{v3}
\ncline[linewidth=0.6mm,linecolor=red]{u4}{v4}
\ncline[linewidth=0.6mm,linecolor=red]{u5}{v5}
\ncline[linewidth=0.6mm,linecolor=red]{u6}{v6}
\ncline[linewidth=0.6mm,linecolor=red]{u7}{v7}

\ncline{u1}{v2}\ncline{u1}{v3}\ncline{u1}{v5}
\ncline{u2}{v3}\ncline{u2}{v4}\ncline{u2}{v6}
\ncline{u3}{v4}\ncline{u3}{v5}\ncline{u3}{v7}
\ncline{u4}{v5}\ncline{u4}{v6}\ncline{u4}{v1}
\ncline{u5}{v6}\ncline{u5}{v7}\ncline{u5}{v2}
\ncline{u6}{v7}\ncline{u6}{v1}\ncline{u6}{v3}
\ncline{u7}{v1}\ncline{u7}{v2}\ncline{u7}{v4}

\end{pspicture}
\caption{{\small The graphs $T_{2n}$ in Lemma \ref{lem-T2n} and $S_{14}$ in Lemma \ref{lem-S14}.}\label{fig1}}
\end{center}
\end{figure}

\begin{lem}{(McKee and Smyth \cite{McSm07})}\label{lem-T2n}
Let $P$ be a permutation matrix of order $n$ such that $P+P^T$
is the adjacency matrix of the cycle $C_n$ and
$$
A_{n}=\left[\begin{array}{cc}
P+P^T & P-P^T \\
P^T-P & -(P+P^T) \\
\end{array}\right].
$$
Then $A_n$ is the adjacency matrix of the $2n$-vertex toroidal tessellation $T_{2n}$ (see Figure \ref{fig1}),
whose eigenvalues are $\pm 2$, each with multiplicity $n$.
\end{lem}

\begin{lem}{(McKee and Smyth \cite{McSm07})}\label{lem-S14}
Let $W(7,4)=(w_{ij})$ be the weighing matrix of order $7$ and weight $4$, where
$w_{ij} = w_{1,\ell}$ for $\ell\equiv j-i+1({\rm mod}\ 7)$ and
$
(w_{11},w_{12},w_{13},w_{14},w_{15},w_{16},w_{17}) = (-1,1,1,0,1,0,0).
$
Let
$$
W(14,4)=
\left[\begin{matrix}
0 & 1 \\
1 & 0 \\
\end{matrix}\right]\otimes W(7,4).
$$
Then $W(14,4)$ is the adjacency matrix of the $14$-vertex signed graph $S_{14}$ (see Figure \ref{fig1})
and its eigenvalues are $\pm 2$ with the same multiplicity $7$.
\end{lem}

\begin{exa}{(Stinson \cite{Stin04})}\label{exa-Kn}
{\rm
For each $n\in \{2,6,10,14,18,26,30\}$, there exists a symmetric conference matrices $W(n,n-1)$.
Then, $W(n,n-1)$ is a signed adjacency matrix of $K_n$ and its eigenvalues are $\pm\sqrt{n-1}$, each with multiplicity $n/2$.}
\end{exa}

By the property of Kronecker product of matrices, we have the following examples.

\begin{exa}\label{example-complete-k-part-graph}
{\rm
Let $W(k,k-1)$ be a symmetric conference matrix of order $k$ and
$H_n$ be a symmetric Hadamard matrix of order $n$. Then $W(k,k-1)\otimes H_n$
is a signed adjacency matrix of the complete $k$-partite graph
$K_{n,n,\dots,n}$ and its eigenvalues are $\pm\sqrt{(k-1)n}$ with the same multiplicity.
In particular, $W(6,5)\otimes H_2$ is a signed adjacency matrix of the complete $6$-partite graph
$K_{2,2,2,2,2,2}$. }
\end{exa}

\begin{exa}\label{example-Hn-AGamma}
{\rm
Let $\Gamma$ be a signed graph of order $m$ and $H_n$ be a symmetric Hadamard matrix of order $n$.
Then $H_n\otimes A(\Gamma)$ is an adjacency matrix of the signed
graph $\Gamma^{(n)}$ of order $mn$ obtained from $\Gamma$.
If $\Gamma$ has exactly two distinct eigenvalues $\pm\theta$,
then $\Gamma^{(n)}$ has exactly two distinct eigenvalues $\pm\theta\sqrt{n}$.
}
\end{exa}

\section{Eigenvalues of signed Cartesian product and signed semi-strong product graphs}

In this section, we discuss the adjacency eigenvalues of
$\Gamma_1\widetilde{\Box} \Gamma_2$ and $\Gamma_1\widetilde{\bowtie} \Gamma_2$
and obtain a sufficient and necessary condition such that the spectrums of
$\Gamma_1\widetilde{\Box}\Gamma_2$ and $\Gamma_1\widetilde{\bowtie}\Gamma_2$
are symmetric.

\begin{thm}\label{thm-eigenvalues}
Let $\Gamma_1=(G_1,\sigma_1)$ be a signed bipartite graph of order $n$ with bipartition $(V_1,V_2)$
and $\Gamma_2=(G_2,\sigma_2)$ be a signed graph of order $m$.
If $\lambda^2$ is an eigenvalue of $A(\Gamma_1)^2$ with multiplicity $p$
and $\mu^2$ is an eigenvalue of $A(\Gamma_2)^2$ with multiplicity $q$, then
each of the following holds.

(i) $\lambda^2 + \mu^2$ (resp. $(\lambda^2 + 1)\mu^2$) is an eigenvalue of
$A(\Gamma_1\widetilde{\Box} \Gamma_2)^2$ (resp. $A(\Gamma_1\widetilde{\bowtie} \Gamma_2)^2$)
with multiplicity $pq$.

(ii) If $\lambda=0$ and $\mu\neq 0$, then $\pm\mu$ are eigenvalues of
$\Gamma_1\widetilde{\Box} \Gamma_2$ (also $\Gamma_1\widetilde{\bowtie} \Gamma_2$)
with multiplicities $\frac{1}{2}pq \pm \frac{1}{2}(n-2|V_1|)(q-2t)$ respectively,
where $t$ is the multiplicity of eigenvalue $\mu$ of $A(\Gamma_2)$.

(iii) If $\lambda\neq 0$, then $\pm \sqrt{\lambda^2 + \mu^2}$ are  eigenvalues of
$\Gamma_1\widetilde{\Box} \Gamma_2$, each with multiplicity $pq/2$.

(iv) If $\lambda \mu\neq 0$, then $\pm \sqrt{(\lambda^2 + 1)\mu^2}$ are eigenvalues of
$\Gamma_1\widetilde{\bowtie} \Gamma_2$, each with multiplicity $pq/2$.
\end{thm}

\begin{cor}\label{cor-eig-G1G2}
For $i=1,2$, let $\Gamma_i$ be a signed graph with exactly two distinct eigenvalues $\pm \theta_i$,
where $\Gamma_1$ is bipartite. Then $\Gamma_1\widetilde{\Box}\Gamma_2$ and $\Gamma_1\widetilde{\bowtie}\Gamma_2$
have exactly two distinct eigenvalues $\pm \sqrt{\theta_1^2+\theta_2^2}$ and $\pm \sqrt{(\theta_1^2+1)\theta_2^2}$, respectively.
\end{cor}

The following theorem gives a
sufficient and necessary condition such that the spectrums of
$\Gamma_1\widetilde{\Box}\Gamma_2$ and $\Gamma_1\widetilde{\bowtie}\Gamma_2$
are symmetric.

\begin{thm}\label{thm-symmetric}
Let $\Gamma_1$ be a signed bipartite graph and $\Gamma_2$ be a signed graph.
The spectrum of $\Gamma_1\widetilde{\Box}\Gamma_2$ (resp. $\Gamma_1\widetilde{\bowtie}\Gamma_2$)
is symmetric if and only if $\Gamma_1$ is balanced or the spectrum of $\Gamma_2$ is symmetric.
\end{thm}

In the following proofs of Theorem \ref{thm-eigenvalues} and Theorem \ref{thm-symmetric}, we always
assume that $A_1:=A(\Gamma_1)=\left[
\begin{matrix}
O_s & P \\
P^T & O_{n-s}
\end{matrix}
\right]$ and $A_2:=A(\Gamma_2)$, where $|V_1|=s$ and $P$ is an $s\times (n-s)$ matrix.

\bigskip
\noindent{\bf Proof of Theorem \ref{thm-eigenvalues} (i).}
By (\ref{e1.1}), we have
\begin{align}\label{e4.1}
A(\Gamma_1\widetilde{\Box} \Gamma_2)^2
= &\left(A_1 \otimes I_m +
\left[
\begin{matrix}
I_s & O \\
O & -I_{n-s}
\end{matrix}
\right] \otimes A_2\right)^2 \nonumber\\ 
= &  A_1^2 \otimes I_m +
\left[
\begin{matrix}
I_s & O \\
O & I_{n-s}
\end{matrix}
\right] \otimes A_2^2 \nonumber\\
 & + \left(\left[
\begin{matrix}
O & P \\
P^T & O
\end{matrix}
\right]\otimes I_m \right) \left( \left[
\begin{matrix}
I_s & O \\
O & -I_{n-s}
\end{matrix}
\right] \otimes A_2 \right) \nonumber\\
& +  \left( \left[
\begin{matrix}
I_s & O \\
O & -I_{n-s}
\end{matrix}
\right] \otimes A_2 \right)
\left(\left[
\begin{matrix}
O & P \\
P^T & O
\end{matrix}
\right]\otimes I_m \right) \nonumber\\
= & A_1^2 \otimes I_m + I_n \otimes A_2^2
 + \left[
\begin{matrix}
O & -P \\
P^T & O
\end{matrix}
\right]\otimes A_2
 +   \left[
\begin{matrix}
O & P \\
-P^T & O
\end{matrix}
\right] \otimes A_2  \nonumber\\ 
= & A_1^2 \otimes I_{m} + I_{n} \otimes A_2^2.
\end{align}
For each $i=1,\dots,p$ and $j=1,\dots,q$, let $X_i$ and $Y_j$ be eigenvectors of $A_1^2$ and $A_2^2$
with respect to eigenvalues $\lambda^2$ and $\mu^2$, respectively.
Thus, by (\ref{e4.1}), we have
$$
A(\Gamma_1\widetilde{\Box} \Gamma_2)^2(X_i\otimes Y_j) = (A_1^2\otimes I_m+I_n\otimes A_2^2)(X_i\otimes Y_j)=(\lambda^2 + \mu^2) (X_i\otimes Y_j).
$$
Therefore, $\lambda^2 + \mu^2$ is an eigenvalue of
$A(\Gamma_1\widetilde{\Box} \Gamma_2)^2$ with multiplicity $pq$.

By (\ref{e1.2}), we have
\begin{align}\label{e4.2}
A(\Gamma_1\widetilde{\bowtie} \Gamma_2)^2
= &\left(\left[
\begin{matrix}
I_s & P \\
P^T & -I_{n-s}
\end{matrix}
\right] \otimes A_2\right)^2 \nonumber\\ 
= & \left(\left[
\begin{matrix}
O_s & P \\
P^T & O_{n-s}
\end{matrix}
\right]+\left[
\begin{matrix}
I_s & O \\
O & -I_{n-s}
\end{matrix}
\right]\right)^2 \otimes A_2^2 \nonumber\\
=& (A_1^2+I_n) \otimes A_2^2 + \left[
\begin{matrix}
O_s & -P+P \\
P^T-P^T & O_{n-s}
\end{matrix}
\right]\otimes A_2^2 \nonumber\\
=& (A_1^2+I_n) \otimes A_2^2.
\end{align}
For each $i=1,\dots,p$ and $j=1,\dots,q$, let $X_i$ and $Y_j$ be eigenvectors of $A_1^2$ and $A_2^2$
with respect to eigenvalues $\lambda^2$ and $\mu^2$, respectively.
Thus, by (\ref{e4.2}), we have
$$
A(\Gamma_1\widetilde{\bowtie} \Gamma_2)^2(X_i\otimes Y_j) =
[(A_1^2 +I_n)\otimes A_2^2](X_i\otimes Y_j)=(\lambda^2 +1) \mu^2 (X_i\otimes Y_j).
$$
Therefore, $(\lambda^2 +1) \mu^2$ is an eigenvalue of
$A(\Gamma_1\widetilde{\bowtie} \Gamma_2)^2$ with multiplicity $pq$. \hfill\rule{2mm}{2mm}

\begin{lem}\label{lem-null-space}
Let $\Gamma$ be a signed bipartite graph of order $n$ with bipartition $(V_1,V_2)$, where $|V_1|=s$, and
$A=\left[
  \begin{matrix}
    O_s & P\\
    P^T & O_{n-s} \\
  \end{matrix}
\right]$ be the adjacency matrix of $\Gamma$.
Let $\{\mathbf{w}_1,\dots,\mathbf{w}_a\}$ be a basis of null space of $P^T$ and
$\{\mathbf{z}_1,\dots,\mathbf{z}_b\}$ be a basis of null space of $P$.
The following $a+b$ vectors of length $n$
$$
\left\{
\left[
  \begin{matrix}
    \mathbf{w}_1 \\
    \bm{0} \\
  \end{matrix}
\right], \dots,
\left[
  \begin{matrix}
    \mathbf{w}_a \\
    \bm{0} \\
  \end{matrix}
\right], \left[
  \begin{matrix}
    \bm{0} \\
    \mathbf{z}_1 \\
  \end{matrix}
\right], \dots,
\left[
  \begin{matrix}
    \bm{0} \\
    \mathbf{z}_b \\
  \end{matrix}
\right]
\right\},
$$
is a basis of null space of $A$.
\end{lem}

\begin{pf}
Since ${\rm rank}(A)={\rm rank}(P)+{\rm rank}(P^T)$, by Rank-Nullity Theorem
$$
n-{\rm rank}(A)=(n-s-{\rm rank}(P))+(s-{\rm rank}(P^T))=a+b.
$$
The result follows.
\end{pf}

\bigskip
\noindent{\bf Proof of Theorem \ref{thm-eigenvalues} (ii).}
Since the multiplicity of eigenvalue $\mu$ of $A_2$ is $t$,
the multiplicity of eigenvalue $-\mu$ of $A_2$ is $q-t$.
Assume that $A_2Y_j = \mu Y_j$ for each $1\leq j\leq t$ and $A_2Y'_k = -\mu Y'_k$ for each $1\leq k\leq q-t$.
In particular, if $t=0$, then $1\leq k\leq q$ and there exists no such $Y_j$;
if $t=q$, then $1\leq j\leq q$ and there exists no such $Y'_k$.

By the assumption, $\lambda=0$ is an eigenvalue of $A_1^2$ (and so $A_1$) with multiplicity $p$.
Hence, the rank of $A_1$ is ${\rm rank}(A_1)=n-p$ and
${\rm rank}(P)={\rm rank}(P^T)=(n-p)/2$. Thus,
the nullity of $P^T$ is
$$
r:= s - {\rm rank}(P^T)= p/2 - (n-2s)/2
$$
and the nullity of $P$ is $p-r = p/2 + (n-2s)/2$.
Suppose that $\{X_{11},\dots,X_{r1}\}$ is a basis of null space of $P^T$ and
$\{X_{12},\dots,X_{(p-r)2}\}$ is a basis of null space of $P$.
Let
$$
Z_i:=\left[
  \begin{matrix}
    X_{i1} \\
    \bm{0} \\
  \end{matrix}
\right],
Z'_{\ell}:=\left[
  \begin{matrix}
    \bm{0} \\
X_{\ell 2} \\
  \end{matrix}
\right]
$$
be column vectors of length $n$ for each $1\leq i\leq r$ and $1\leq \ell\leq p-r$.
In particular, if $r=0$, then $1\leq \ell\leq p=n-2s$ and there is no such $Z_i$;
if $r=p$, then $1\leq i\leq p=2s-n$ and there is no such $Z'_{\ell}$.
By Lemma \ref{lem-null-space},
$\{Z_1,\dots,Z_r\}\cup \{Z'_1,\dots,Z'_{p-r}\}$ is a basis of
null space of $A_1$. Therefore, $A_1Z_i=A_1Z'_{\ell}=\bm{0}$ and
for every $1\leq i\leq r$, $1\leq j\leq t$ and $1\leq k\leq q-t$,
\begin{align*}
A(\Gamma_1\widetilde{\Box} \Gamma_2)(Z_i\otimes Y_j)&=\mu (Z_i\otimes Y_j)=A(\Gamma_1\widetilde{\bowtie} \Gamma_2)(Z_i\otimes Y_j),\\
A(\Gamma_1\widetilde{\Box} \Gamma_2)(Z_i\otimes Y'_k)&=-\mu (Z_i\otimes Y'_k)=A(\Gamma_1\widetilde{\bowtie} \Gamma_2)(Z_i\otimes Y'_k).
\end{align*}
For every $1\leq \ell\leq p-r$, $1\leq j\leq t$ and $1\leq k\leq q-t$,
\begin{align*}
A(\Gamma_1\widetilde{\Box} \Gamma_2)(Z'_{\ell}\otimes Y'_k)&=\mu (Z'_{\ell}\otimes Y'_k)
 =A(\Gamma_1\widetilde{\bowtie} \Gamma_2)(Z'_{\ell}\otimes Y'_k),\\
A(\Gamma_1\widetilde{\Box} \Gamma_2)(Z'_{\ell}\otimes Y_j)&=-\mu (Z'_{\ell}\otimes Y_j)
 =A(\Gamma_1\widetilde{\bowtie} \Gamma_2)(Z'_{\ell}\otimes Y_j).
\end{align*}
Note that all of $Z_i$, $Z'_{\ell}$, $Y_j$ and $Y'_k$ are nonzero vectors
for each $1\leq i\leq r$, $1\leq \ell\leq p-r$, $1\leq j\leq t$ and $1\leq k\leq q-t$.
Hence, the Kronecker products of them are also nonzero vectors.
By $(Z_i\otimes Y_j)^T (Z'_{\ell}\otimes Y'_k)=0$, we have $Z_i\otimes Y_j$
and $Z'_{\ell}\otimes Y'_k$ are
$$
rt+(p-r)(q-t) =pq/2 + (n-2s)(q-2t)/2
$$
eigenvectors of $A(\Gamma_1\widetilde{\Box} \Gamma_2)$ (resp. $A(\Gamma_1\widetilde{\bowtie} \Gamma_2)$)
with respect to eigenvalue $\mu$. By $(Z_i\otimes Y'_k)^T (Z'_{\ell}\otimes Y_j)=0$,
we know that $Z_i\otimes Y'_k$ and $Z'_{\ell}\otimes Y_j$ are
$$
r(q-t)+(p-r)t = pq/2 - (n-2s)(q-2t)/2
$$
eigenvectors of $A(\Gamma_1\widetilde{\Box} \Gamma_2)$
(resp. $A(\Gamma_1\widetilde{\bowtie} \Gamma_2)$) with respect to eigenvalue $-\mu$.
Thus, $\pm\mu$ are eigenvalues of
$\Gamma_1\widetilde{\Box} \Gamma_2$ (resp. $\Gamma_1\widetilde{\bowtie} \Gamma_2$)
with multiplicities $\frac{1}{2}pq\pm \frac{1}{2}(n-2s)(q-2t)$, respectively. \hfill\rule{2mm}{2mm}

\bigskip
\noindent{\bf Proof of Theorem \ref{thm-eigenvalues} (iii).}
Suppose that $\lambda\neq 0$.
Since $\Gamma_1$ is bipartite, $\lambda$ and $-\lambda$ are
eigenvalues of $\Gamma_1$, each with multiplicity $p/2$.
Without loss of generality, assume that $\mu\geq 0$,
$A_2Y_j = \mu Y_j$ for each $j=1,\dots,t$ and $A_2Y'_k = -\mu Y'_k$ for each $k=1,\dots,q-t$.
In particular, if $t=0$, then $1\leq k\leq q$ and there exists no such $Y_j$;
if $t=q$, then $1\leq j\leq q$ and there exists no such $Y'_k$.
Note that if $\mu=0$, then $t=q$.
Now, for $i=1,\dots,p/2$, suppose that
$
X_i=\left[\begin{matrix}
X_{i1} \\
X_{i2}
\end{matrix}\right]
$
is the unit vector such that $A_1X_i=\lambda X_i$, where $X_{i1}$ and $X_{i2}$
are column vectors of length $s$ and $n-s$ respectively.
Then $PX_{i2}=\lambda X_{i1}$ and $P^TX_{i1}=\lambda X_{i2}$.
For each $i=1,\dots,p/2$, let
$X'_i=\left[\begin{matrix}
X_{i1} \\
-X_{i2}
\end{matrix}\right]
$, then $A_1X'_i=-\lambda X'_i$. Since $\lambda\neq 0$, we have
$X_i^TX'_i=0$ and so $X_{i1}^TX_{i1}=X_{i2}^TX_{i2}=\frac{1}{2}$,
which implies that $X_{i1}$ and $X_{i2}$ are nonzero vectors.
Based on eigenvalues $\pm\lambda, \pm\mu$ and the corresponding eigenvectors, we construct $pq/2$ vectors as follows
$$
Z_i\otimes Y_j=\left[
\begin{matrix}
(\sqrt{\lambda^2+\mu^2}+\mu) X_{i1} \\
\lambda X_{i2} \\
\end{matrix}
\right]\otimes Y_j,  \, \,
W_i\otimes Y'_k=\left[
\begin{matrix}
\lambda  X_{i1} \\
(\sqrt{\lambda^2+\mu^2}+\mu) X_{i2} \\
\end{matrix}
\right]\otimes Y'_k,
$$
for each $i=1,\dots,p/2$, $j=1,\dots,t$ and $k=1,\dots,q-t$, and construct $pq/2$ vectors as follows
$$
Z'_i\otimes Y_j=\left[
\begin{matrix}
-\lambda X_{i1} \\
(\sqrt{\lambda^2+\mu^2}+\mu)X_{i2} \\
\end{matrix}
\right]\otimes Y_j,  \, \,
W'_i\otimes Y'_k=\left[
\begin{matrix}
(\sqrt{\lambda^2+\mu^2}+\mu) X_{i1} \\
-\lambda X_{i2} \\
\end{matrix}
\right]\otimes Y'_k,
$$
for $i=1,\dots,p/2$, $j=1,\dots,t$ and $k=1,\dots,q-t$.
Then, we have
\begin{align}\label{e4.3}
A(\Gamma_1\widetilde{\Box} \Gamma_2)\cdot(Z_i\otimes Y_j) = &
\left[
\begin{matrix}
I_s\otimes A_2 & P\otimes I_m \\
P^T\otimes I_m & -I_{n-s}\otimes A_2
\end{matrix}
\right] \cdot \left[
\begin{matrix}
(\sqrt{\lambda^2+\mu^2}+\mu) X_{i1} \otimes Y_j \\
\lambda X_{i2} \otimes Y_j
\end{matrix}
\right]  \nonumber\\ 
= & \left[
\begin{matrix}
(\sqrt{\lambda^2+\mu^2}+\mu)X_{i1}\otimes \mu Y_j +\lambda^2 X_{i1}\otimes Y_j \\
(\sqrt{\lambda^2+\mu^2}+\mu)\lambda X_{i2}\otimes Y_j -\lambda X_{i2}\otimes \mu Y_j \\
\end{matrix}
\right] \nonumber\\ 
= & \sqrt{\lambda^2+\mu^2}\cdot \left[
\begin{matrix}
(\sqrt{\lambda^2+\mu^2}+\mu) X_{i1}\otimes  Y_j \\
\lambda X_{i2}\otimes  Y_j \\
\end{matrix}
\right] \nonumber\\ 
= & \sqrt{\lambda^2+\mu^2}\cdot (Z_i\otimes Y_j),
\end{align}
\begin{align}\label{e4.4}
A(\Gamma_1\widetilde{\Box} \Gamma_2)\cdot(W_i\otimes Y'_k) = &
\left[
\begin{matrix}
I_s\otimes A_2 & P\otimes I_m \\
P^T\otimes I_m & -I_{n-s}\otimes A_2
\end{matrix}
\right]\cdot \left[
\begin{matrix}
\lambda  X_{i1}\otimes Y'_k \\
(\sqrt{\lambda^2+\mu^2}+\mu) X_{i2}\otimes Y'_k \\
\end{matrix}
\right] \nonumber\\ 
= & \left[
\begin{matrix}
-\lambda X_{i1}\otimes \mu Y'_k + (\sqrt{\lambda^2+\mu^2}+\mu)\lambda X_{i1}\otimes Y'_k\\
\lambda^2 X_{i2}\otimes Y'_k + (\sqrt{\lambda^2+\mu^2}+\mu) X_{i2}\otimes \mu Y'_k
\end{matrix}
\right]\nonumber\\ 
= & \sqrt{\lambda^2+\mu^2}\cdot \left[
\begin{matrix}
\lambda X_{i1}\otimes Y'_k \\
(\sqrt{\lambda^2+\mu^2}+\mu) X_{i2}\otimes Y'_k \\
\end{matrix}
\right] \nonumber\\ 
= & \sqrt{\lambda^2+\mu^2}\cdot (W_i\otimes Y'_k),
\end{align}
\begin{align}\label{e4.5}
A(\Gamma_1\widetilde{\Box} \Gamma_2)\cdot(Z'_i\otimes Y_j) = &
\left[
\begin{matrix}
I_s\otimes A_2 & P\otimes I_m \\
P^T\otimes I_m & -I_{n-s}\otimes A_2
\end{matrix}
\right] \cdot \left[
\begin{matrix}
-\lambda X_{i1}\otimes Y_j \\
(\sqrt{\lambda^2+\mu^2}+\mu)X_{i2}\otimes Y_j \\
\end{matrix}
\right] \nonumber\\ 
= & \left[
\begin{matrix}
-\lambda X_{i1}\otimes \mu Y_j + (\sqrt{\lambda^2+\mu^2}+\mu)\lambda X_{i1}\otimes Y_j\\
-\lambda^2 X_{i2}\otimes Y_j - (\sqrt{\lambda^2+\mu^2}+\mu) X_{i2}\otimes \mu Y_j
\end{matrix}
\right]\nonumber\\ 
= & -\sqrt{\lambda^2+\mu^2}\cdot \left[
\begin{matrix}
-\lambda X_{i1} \otimes  Y_j \\
(\sqrt{\lambda^2+\mu^2}+\mu)X_{i2} \otimes  Y_j \\
\end{matrix}
\right] \nonumber\\
= & -\sqrt{\lambda^2+\mu^2}\cdot (Z'_i\otimes Y_j),
\end{align}
\begin{align}\label{e4.6}
A(\Gamma_1\widetilde{\Box} \Gamma_2)\cdot(W'_i\otimes Y'_k) = &
\left[
\begin{matrix}
I_s\otimes A_2 & P\otimes I_m \\
P^T\otimes I_m & -I_{n-s}\otimes A_2
\end{matrix}
\right] \cdot \left[
\begin{matrix}
(\sqrt{\lambda^2+\mu^2}+\mu) X_{i1}\otimes Y'_k \\
-\lambda X_{i2}\otimes Y'_k \\
\end{matrix}
\right] \nonumber\\ 
= & \left[
\begin{matrix}
-(\sqrt{\lambda^2+\mu^2}+\mu) X_{i1} \otimes \mu Y'_k -\lambda^2 X_{i1}\otimes Y'_k\\
(\sqrt{\lambda^2+\mu^2}+\mu) \lambda X_{i2}\otimes Y'_k -\lambda X_{i2}\otimes \mu Y'_k
\end{matrix}
\right] \nonumber\\ 
= & -\sqrt{\lambda^2+\mu^2}\cdot \left[
\begin{matrix}
(\sqrt{\lambda^2+\mu^2}+\mu) X_{i1}\otimes  Y'_k \\
-\lambda X_{i2}\otimes  Y'_k \\
\end{matrix}
\right] \nonumber\\
= & -\sqrt{\lambda^2+\mu^2}\cdot (W'_i\otimes Y'_k).
\end{align}
Since $\lambda\neq 0$, and $X_{i1}$ and $X_{i2}$ are nonzero, we know that
all of $Z_i, W_i, Z'_i, W'_i, Y_j, Y'_k$ are nonzero vectors
for each $i\in \{1,\dots,p/2\}$, $j\in \{1,\dots,t\}$ and $k\in \{1,\dots,q-t\}$,
and the Kronecker products of them are also nonzero vectors.
As $(Z_i\otimes Y_j)^T (W_i\otimes Y'_k)=0$, by (\ref{e4.3}) and (\ref{e4.4}), we have $Z_i\otimes Y_j$
and $W_i\otimes Y'_k$ are $pq/2$ eigenvectors
of $A(\Gamma_1\widetilde{\Box} \Gamma_2)$ with respect to eigenvalue $\sqrt{\lambda^2+\mu^2}$. As
$(Z'_i\otimes Y_j)^T (W'_i\otimes Y'_k)=0$, by (\ref{e4.5}) and (\ref{e4.6}), we have $Z'_i\otimes Y_j$
and $W'_i\otimes Y'_k$ are $pq/2$ eigenvectors
of $A(\Gamma_1\widetilde{\Box} \Gamma_2)$ with respect to eigenvalue $-\sqrt{\lambda^2+\mu^2}$.
Therefore, $\pm\sqrt{\lambda^2+\mu^2}$ are adjacency eigenvalues of $\Gamma_1\widetilde{\Box} \Gamma_2$,
each with multiplicity $pq/2$. \hfill\rule{2mm}{2mm}

\bigskip
\noindent{\bf Proof of Theorem \ref{thm-eigenvalues} (iv).}
Suppose that $\lambda\mu\neq 0$.
Since $\Gamma_1$ is bipartite, $\lambda$ and $-\lambda$ are
eigenvalues of $\Gamma_1$, each with multiplicity $p/2$.
Without loss of generality, assume that
$A_2Y_j = \mu Y_j$ for each $j=1,\dots,t$ and $A_2Y'_k = -\mu Y'_k$ for each $k=1,\dots,q-t$.
In particular, if $t=0$, then $1\leq k\leq q$ and there exists no such $Y_j$;
if $t=q$, then $1\leq j\leq q$ and there exists no such $Y'_k$.
Now, for each $i=1,\dots,p/2$, suppose that
$
X_i=\left[\begin{matrix}
X_{i1} \\
X_{i2}
\end{matrix}\right]
$
is the unit vector such that $A_1X_i=\lambda X_i$, where $X_{i1}$ and $X_{i2}$
are column vectors of length $s$ and $n-s$ respectively.
Then $PX_{i2}=\lambda X_{i1}$ and $P^TX_{i1}=\lambda X_{i2}$.
Let $
X'_i=\left[\begin{matrix}
X_{i1} \\
-X_{i2}
\end{matrix}\right]
$, then $A_1X'_i=-\lambda X'_i$ and so $X_i^TX'_i=0$. Thus
$
(X_{i1})^TX_{i1}=(X_{i2})^TX_{i2}=\frac{1}{2},
$
and so $X_{i1}$ and $X_{i2}$ are nonzero vectors.
Based on eigenvalues $\pm\lambda, \pm\mu$ and the corresponding eigenvectors, we construct $pq/2$ vectors as follows
$$
Z_i\otimes Y_j=\left[
\begin{matrix}
(\sqrt{\lambda^2+1}+1) X_{i1} \\
\lambda X_{i2} \\
\end{matrix}
\right]\otimes Y_j,  \, \,
Z'_i\otimes Y'_k=\left[
\begin{matrix}
-\lambda  X_{i1} \\
(\sqrt{\lambda^2+1}+1) X_{i2} \\
\end{matrix}
\right]\otimes Y'_k,
$$
for each $i=1,\dots,p/2$, $j=1,\dots,t$ and $k=1,\dots,q-t$, and construct $pq/2$ vectors as follows
$$
Z'_i\otimes Y_j=\left[
\begin{matrix}
-\lambda X_{i1} \\
(\sqrt{\lambda^2+1}+1)X_{i2} \\
\end{matrix}
\right]\otimes Y_j,  \, \,
Z_i\otimes Y'_k=\left[
\begin{matrix}
(\sqrt{\lambda^2+1}+1) X_{i1} \\
\lambda X_{i2} \\
\end{matrix}
\right]\otimes Y'_k,
$$
for each $i=1,\dots,p/2$, $j=1,\dots,t$ and $k=1,\dots,q-t$.
Since
\begin{align*}
\left[
\begin{matrix}
I_s & P \\
P^T & -I_{n-s}
\end{matrix}
\right]Z_i  = &
\left[
\begin{matrix}
I_s & P \\
P^T & -I_{n-s}
\end{matrix}
\right]  \left[
\begin{matrix}
(\sqrt{\lambda^2+1}+1) X_{i1} \\
\lambda X_{i2} \\
\end{matrix}
\right]  \\ 
= & \left[
\begin{matrix}
(\sqrt{\lambda^2+1}+1) X_{i1} + \lambda^2 X_{i1} \\
(\sqrt{\lambda^2+1}+1) \lambda X_{i2}-\lambda X_{i2} \\
\end{matrix}
\right]  \\ 
= & \sqrt{\lambda^2+1}Z_i , \\
\left[
\begin{matrix}
I_s & P \\
P^T & -I_{n-s}
\end{matrix}
\right]Z'_i = &
\left[
\begin{matrix}
I_s & P \\
P^T & -I_{n-s}
\end{matrix}
\right] \left[
\begin{matrix}
-\lambda X_{i1} \\
(\sqrt{\lambda^2+1}+1)X_{i2} \\
\end{matrix}
\right]  \\ 
= & \left[
\begin{matrix}
-\lambda X_{i1} + (\sqrt{\lambda^2+1}+1) \lambda X_{i1} \\
- \lambda^2 X_{i2} - (\sqrt{\lambda^2+1}+1) X_{i2} \\
\end{matrix}
\right]  \\ 
= & -\sqrt{\lambda^2+1} Z'_i ,
\end{align*}
we can obtain the following equations
\begin{align}
A(\Gamma_1\widetilde{\bowtie} \Gamma_2)\cdot(Z_i\otimes Y_j) &= \sqrt{\lambda^2+1}Z_i\otimes A_2 Y_j
=  \mu\sqrt{\lambda^2+1}\cdot (Z_i\otimes Y_j), \label{e4.7}\\
A(\Gamma_1\widetilde{\bowtie} \Gamma_2)\cdot(Z'_i\otimes Y'_k) &= -\sqrt{\lambda^2+1}Z'_i\otimes A_2 Y'_k
=  \mu\sqrt{\lambda^2+1}\cdot (Z'_i\otimes Y'_k), \label{e4.8}\\
A(\Gamma_1\widetilde{\bowtie} \Gamma_2)\cdot(Z'_i\otimes Y_j) &= -\sqrt{\lambda^2+1}Z'_i \otimes A_2 Y_j
=  -\mu\sqrt{\lambda^2+1}\cdot (Z'_i\otimes Y_j), \label{e4.9}\\
A(\Gamma_1\widetilde{\bowtie} \Gamma_2)\cdot(Z_i\otimes Y'_k) &= \sqrt{\lambda^2+1}Z_i \otimes A_2 Y'_k
=  -\mu\sqrt{\lambda^2+1}\cdot (Z_i\otimes Y'_k).\label{e4.10}
\end{align}
Since $\lambda\mu\neq 0$, $X_{i1}$ and $X_{i2}$ are nonzero, we know that
all of $Z_i, Z'_i, Y_j, Y'_k$ are nonzero vectors
for each $i\in \{1,\dots,p/2\}$, $j\in \{1,\dots,t\}$ and $k\in \{1,\dots,q-t\}$, and
the Kronecker products of them are also nonzero.
As $(Z_i\otimes Y_j)^T (Z'_i\otimes Y'_k)=0$, by (\ref{e4.7}) and (\ref{e4.8}), we have $Z_i\otimes Y_j$
and $Z'_i\otimes Y'_k$ are $pq/2$ eigenvectors
of $A(\Gamma_1\widetilde{\bowtie} \Gamma_2)$ with respect to eigenvalue $\mu\sqrt{\lambda^2+1}$. As
$(Z'_i\otimes Y_j)^T (Z_i\otimes Y'_k)=0$, by (\ref{e4.9}) and (\ref{e4.10}), we have $Z'_i\otimes Y_j$
and $Z_i\otimes Y'_k$ are $pq/2$ eigenvectors
of $A(\Gamma_1\widetilde{\bowtie} \Gamma_2)$ with respect to eigenvalue $-\mu\sqrt{\lambda^2+1}$.
Therefore, $\pm\mu\sqrt{\lambda^2+1}$ are adjacency eigenvalues of $\Gamma_1\widetilde{\bowtie} \Gamma_2$,
each with multiplicity $pq/2$.  \hfill\rule{2mm}{2mm}

\bigskip
\noindent{\bf Proof of Theorem \ref{thm-symmetric}.}
Let $\lambda^2$ be any eigenvalue of $A(\Gamma_1)^2$ with multiplicity $p$
and $\mu^2$ be any eigenvalue of $A(\Gamma_2)^2$ with multiplicity $q$,
where $\mu$ is the eigenvalue of $\Gamma_2$ with multiplicity $t$.

(a) Consider $\Gamma_1\widetilde{\Box} \Gamma_2$.
By Theorem \ref{thm-eigenvalues} (i), $\lambda^2 + \mu^2$ is an eigenvalue of
$A(\Gamma_1\widetilde{\Box} \Gamma_2)^2$ with multiplicity $pq$.

Assume that $\Gamma_1$ is balanced or the spectrum of $\Gamma_2$ is symmetric.
Then the bipartition $(V_1,V_2)$ of $\Gamma_1$ satisfies
$|V_1|= \frac{n}{2}$ or the multiplicity of eigenvalue $\mu (\neq 0)$ of $\Gamma_2$ is equal to $t=\frac{q}{2}$, and so
$
(n-2|V_1|)(q-2t)=0.
$
It suffices to prove that the multiplicities of eigenvalues $\pm\sqrt{\lambda^2+\mu^2}$
of $\Gamma_1\widetilde{\Box} \Gamma_2$ are equal to $\frac{1}{2}pq$ when $\lambda^2+\mu^2\neq 0$.
If $\lambda\neq 0$, then by Theorem \ref{thm-eigenvalues} (iii),
the multiplicities of eigenvalues $\pm\sqrt{\lambda^2+\mu^2}$
of $\Gamma_1\widetilde{\Box} \Gamma_2$ are equal to $\frac{1}{2}pq$.
If $\lambda=0$ and $\mu\neq 0$, then by Theorem \ref{thm-eigenvalues} (ii),
the multiplicities of eigenvalues $\pm\mu$ of $\Gamma_1\widetilde{\Box} \Gamma_2$
are equal to $\frac{1}{2}pq$. Thus, the spectrum of $\Gamma_1\widetilde{\Box} \Gamma_2$ is symmetric.

Conversely, assume that the spectrum of $\Gamma_1\widetilde{\Box} \Gamma_2$ is symmetric.
If all the eigenvalues of $\Gamma_1$ are nonzero, then the rank of
$A(\Gamma_1)$ is $n$ and so ${\rm rank}(P)={\rm rank}(P^T)=\frac{n}{2}$.
This implies $|V_1|=|V_2|$ and so $\Gamma_1$ is balanced.
If $\lambda=0$ and $\mu\neq 0$, then the multiplicities of eigenvalues
$\pm\mu$ of $\Gamma_1\widetilde{\Box} \Gamma_2$ must be equal. By Theorem \ref{thm-eigenvalues} (ii),
we have
$$
pq + (n-2|V_1|)(q-2t) = pq - (n-2|V_1|)(q-2t),
$$
that is $(n-2|V_1|)(q-2t)=0$ and so $\Gamma_1$ is balanced or the spectrum of $\Gamma_2$ is symmetric.

(b) Consider $\Gamma_1\widetilde{\bowtie} \Gamma_2$.
By Theorem \ref{thm-eigenvalues} (i), $(\lambda^2 + 1)\mu^2$ is an eigenvalue of
$A(\Gamma_1\widetilde{\bowtie} \Gamma_2)^2$ with multiplicity $pq$.

Assume that $\Gamma_1$ is balanced or the spectrum of $\Gamma_2$ is symmetric.
Then the bipartition $(V_1,V_2)$ of $\Gamma_1$ satisfies
$|V_1|= \frac{n}{2}$ or the multiplicity of eigenvalue $\mu (\neq 0)$ of $\Gamma_2$ is equal to $t=\frac{q}{2}$, and so
$
(n-2|V_1|)(q-2t)=0.
$
It suffices to prove that the multiplicities of eigenvalues $\pm\sqrt{(\lambda^2+1)\mu^2}$
of $\Gamma_1\widetilde{\bowtie} \Gamma_2$ are equal to $\frac{1}{2}pq$ when $\mu^2\neq 0$.
If $\lambda\neq 0$ and $\mu\neq 0$, then by Theorem \ref{thm-eigenvalues} (iv),
the multiplicities of eigenvalues $\pm\sqrt{(\lambda^2+1)\mu^2}$
of $\Gamma_1\widetilde{\bowtie} \Gamma_2$ are equal to $\frac{1}{2}pq$.
If $\lambda=0$ and $\mu\neq 0$, then by Theorem \ref{thm-eigenvalues} (ii),
the multiplicities of eigenvalues $\pm\mu$ of $\Gamma_1\widetilde{\bowtie} \Gamma_2$
are equal to $\frac{1}{2}pq$. Thus, the spectrum of $\Gamma_1\widetilde{\bowtie} \Gamma_2$ is symmetric.

Conversely, assume that the spectrum of $\Gamma_1\widetilde{\bowtie} \Gamma_2$ is symmetric.
If all the eigenvalues of $\Gamma_1$ are nonzero, then the rank of
$A(\Gamma_1)$ is $n$ and so ${\rm rank}(P)={\rm rank}(P^T)=\frac{n}{2}$.
This implies $|V_1|=|V_2|$ and so $\Gamma_1$ is balanced.
If $\lambda=0$ and $\mu\neq 0$, then the multiplicities of eigenvalues
$\pm\mu$ of $\Gamma_1\widetilde{\bowtie} \Gamma_2$ must be equal. By Theorem \ref{thm-eigenvalues} (ii),
we have
$$
pq + (n-2|V_1|)(q-2t) = pq - (n-2|V_1|)(q-2t),
$$
that is $(n-2|V_1|)(q-2t)=0$ and so $\Gamma_1$ is balanced or the spectrum of $\Gamma_2$ is symmetric.
\hfill\rule{2mm}{2mm}

\section{Induced subgraphs of the signed product graphs}

In this section, we mainly give the proof of Theorem \ref{thm-Delta-2-min} and generalize it to
signed product of $n$ $(n\geq 3)$ graphs.
To establish Theorem \ref{thm-Delta-2-min}, we need the following lemmas.

\begin{lem}(Cauchy's Interlacing Theorem \cite{BrHa12})\label{lem-cauchy}
Let $A$ be an $n\times n$ symmetric matrix, and $B$ be an $m\times m$ principle submatrix of $A$, where $m<n$.
If the eigenvalues of $A$ are $\lambda_1\geq \lambda_2\geq \cdots\geq \lambda_n$, and the eigenvalues of $B$
are $\mu_1\geq \mu_2\geq \cdots\geq \mu_m$, then for all $1\leq i\leq m$,
$$
\lambda_i\geq \mu_i\geq \lambda_{n-m+i}.
$$
\end{lem}

\begin{lem}\label{lem-Delta-sign}
Suppose $\Gamma=(G,\sigma)$ is a signed graph of order $n$, and $A=(a^{\sigma}_{ij})$ is the adjacency matrix of $\Gamma$.
Let $\widetilde{A}=(\tilde{a}_{ij})$ be an $n\times n$ symmetric matrix with
$|\tilde{a}_{ij}|\leq |a^{\sigma}_{ij}|$ for any $1\leq i,j\leq n$. Then
$$
\Delta(\Gamma)\geq \lambda_1(\widetilde{A}).
$$
In particular, $\Delta(\Gamma)\geq \lambda_1(\Gamma)$ when $\widetilde{A}=A$.

\end{lem}

\begin{pf}
It suffices to consider that $\widetilde{A}$ is not an all zero matrix. Thus, $\lambda_1(\widetilde{A})>0$.
Suppose $X=(x_1,x_2,\dots,x_n)^T$ is an eigenvector corresponding to $\lambda_1(\widetilde{A})$.
Then $\lambda_1(\widetilde{A}) X = \widetilde{A} X$. Assume that
$|x_u|=\max\{|x_1|,|x_2|,\dots,|x_n|\}$. Then $|x_u|>0$ and
$$
|\lambda_1(\widetilde{A})x_u|=\left|\sum_{j=1}^n\tilde{a}_{uj}x_j\right| = \left|\sum_{j\sim u}\tilde{a}_{uj}x_j\right|
\leq\sum_{j\sim u}\left|\tilde{a}_{uj}\right||x_u|\leq\sum_{j\sim u}\left|a^{\sigma}_{uj}\right||x_u| \leq \Delta(\Gamma)|x_u|.
$$
Hence, $\Delta(\Gamma)\geq \lambda_1(\widetilde{A})$.
\end{pf}

\begin{lem}\label{lem-Delta-lambda}
Let $\Gamma$ be a connected signed graph of order $n$ with $k$ nonnegative adjacency eigenvalues
$\lambda_1(\Gamma)\geq \cdots\geq \lambda_k(\Gamma)\geq 0$. If $H$ is an $(n-k+1)$-vertex induced subgraph of $\Gamma$, then
$$\Delta(H)\geq \lceil\lambda_k(\Gamma)\rceil.$$
\end{lem}

\begin{pf}
Note that $A(H)$ is an $(n-k+1)\times (n-k+1)$ submatrix of $A(\Gamma)$.
By Lemma \ref{lem-cauchy}, $\lambda_1(H) \geq \lambda_k(\Gamma)$.
By Lemma \ref{lem-Delta-sign}, $\Delta(H)\geq \lambda_1(H) \geq \lambda_k(\Gamma)$.
Hence, $\Delta(H)\geq \lceil\lambda_k(\Gamma)\rceil$.
\end{pf}

\begin{exa}
{\rm
The Petersen graph $(PG,+)$ has spectrum $3^{(1)}, 1^{(5)}, -2^{(4)}$.
If $H$ is a $5$-vertex induced subgraph of $(PG,+)$, then by Lemma \ref{lem-Delta-lambda},
$\Delta(H)\geq 1$ and there exists a subgraph such that the bound is tight.

The signed Petersen graph $(PG,-)$ has spectrum $2^{(4)}, -1^{(5)}, -3^{(1)}$.
If $H$ is a $7$-vertex induced subgraph of $(PG,-)$, then by Lemma \ref{lem-Delta-lambda},
$\Delta(H)\geq 2$ and there exists a subgraph such that the bound is tight.
}
\end{exa}

\noindent{\bf Proof of Theorem \ref{thm-Delta-2-min}.}
Denote $N=mn$.
Let $\Gamma=\Gamma_1\widetilde{\Box}\Gamma_2$ and $\Gamma'=\Gamma_1\widetilde{\bowtie}\Gamma_2$.
Then $H$ (resp. $H'$) is a $(\lfloor\frac{N}{2}\rfloor+1)$-vertex
induced subgraph of $\Gamma$ (resp. $\Gamma'$). By Lemma \ref{lem-Delta-lambda},
$$
\Delta(H) \geq \lambda_{\lceil\frac{N}{2}\rceil}(\Gamma)\  \text{and} \
\Delta(H') \geq \lambda_{\lceil\frac{N}{2}\rceil}(\Gamma').
$$
By Theorem \ref{thm-eigenvalues} (i), $\lambda^2+\mu^2$
is the minimum eigenvalue of $A(\Gamma)^2$ and $(\lambda^2+1)\mu^2$
is the minimum eigenvalue of $A(\Gamma')^2$. Thus, by Theorem \ref{thm-symmetric},
the adjacency spectrums of $\Gamma$ and $\Gamma'$ are symmetric and so
$$
\lambda_{\lceil\frac{N}{2}\rceil}(\Gamma)=\sqrt{\lambda^2+\mu^2}\  \text{and} \
\lambda_{\lceil\frac{N}{2}\rceil}(\Gamma')=\sqrt{(\lambda^2+1)\mu^2}.
$$
Combining these (in)equalities, the results follow. \hfill\rule{2mm}{2mm}

\bigskip
Now, we generalize the signed Cartesian product and signed semi-strong product
of two signed graphs to the product of $n$ signed graphs.

\begin{Def}\label{Def-product-n-graphs}
For $i=1,2,\dots,n-1$, let $\Gamma_i$ be a signed bipartite graph
and $\Gamma_n$ be a signed graph. Let $\Gamma_{\widetilde{\Box},R}^1=\Gamma_{\widetilde{\bowtie},R}^1=\Gamma_1$
and $\Gamma_{\widetilde{\Box},L}^1=\Gamma_{\widetilde{\bowtie},L}^1=\Gamma_n$. For $2\leq k\leq n$,
we define

(i) $\Gamma_{\widetilde{\Box},R}^k=\Gamma_{\widetilde{\Box},R}^{k-1}\widetilde{\Box} \Gamma_k$ and
$\Gamma_{\widetilde{\bowtie},R}^k=\Gamma^{k-1}_{\widetilde{\bowtie},R}\widetilde{\bowtie} \Gamma_k$;

(ii) $\Gamma_{\widetilde{\Box},L}^k=\Gamma_{n-k+1}\widetilde{\Box}\Gamma_{\widetilde{\Box},L}^{k-1}$ and
$\Gamma_{\widetilde{\bowtie},L}^k=\Gamma_{n-k+1}\widetilde{\bowtie}\Gamma_{\widetilde{\bowtie},L}^{k-1}$.
\end{Def}

To illustrate Definition \ref{Def-product-n-graphs}, one can consider $n=3$, that is
\begin{align*}
&\Gamma_{\widetilde{\Box},R}^3= ((\Gamma_1\widetilde{\Box}\Gamma_2)\widetilde{\Box}\Gamma_3) \text{\ and\ }
\Gamma_{\widetilde{\bowtie},R}^3= ((\Gamma_1\widetilde{\bowtie}\Gamma_2)\widetilde{\bowtie}\Gamma_3), \\
& \Gamma_{\widetilde{\Box},L}^3= (\Gamma_1\widetilde{\Box}(\Gamma_2\widetilde{\Box}\Gamma_3)) \text{\ and\ }
\Gamma_{\widetilde{\bowtie},L}^3= (\Gamma_1\widetilde{\bowtie}(\Gamma_2\widetilde{\bowtie}\Gamma_3)).
\end{align*}
By Lemma \ref{lem-HaIK11} and Lemma \ref{lem-GaRW76}, the Cartesian product and semi-strong product
of two bipartite graphs are still bipartite. Therefore, Definition \ref{Def-product-n-graphs} (i)
is well-defined. Since the Kronecker product of matrices is an associative operation, the underling graphs of
$\Gamma_{\widetilde{\Box},R}^n$ and $\Gamma_{\widetilde{\Box},L}^n$ are isomorphic.
However, by Lemma \ref{lem-GaRW76} (iii) and Corollary \ref{cor-GaRW76},
the underling graphs of $\Gamma_{\widetilde{\bowtie},R}^n$ and $\Gamma_{\widetilde{\bowtie},L}^n$ are not isomorphic.

\bigskip
By Definition \ref{Def-product-n-graphs}, Theorem \ref{thm-eigenvalues} (i) can be easily generalized to the following theorem.

\begin{thm}\label{thm-eigenvalues-n-graph}
For $i=1,2,\dots,n$, let $\Gamma_i=(G_i,\sigma_i)$ be a signed graph
and $\theta_i^2$ be an eigenvalue of $A(\Gamma_i)^2$ with multiplicity $p_i$,
where $\Gamma_1,\dots,\Gamma_{n-1}$ are bipartite. Then
$\sum_{i=1}^n\theta_i^2$, $\sum_{i=1}^n\theta_i^2$,
$\theta_n^2\prod_{i=1}^{n-1}(\theta_i^2+1)$
and $\sum_{k=1}^n\prod_{i=k}^{n}\theta_i^2$ are eigenvalues of
$\Gamma_{\widetilde{\Box},L}^n$, $\Gamma_{\widetilde{\Box},R}^n$,
$\Gamma_{\widetilde{\bowtie},L}^n$ and $\Gamma_{\widetilde{\bowtie},R}^n$
with multiplicity $p_1p_2\cdots p_n$, respectively.
\end{thm}

By Theorem \ref{thm-eigenvalues-n-graph}, we have the following corollary immediately.

\begin{cor}\label{cor-eig-G1Gn}
For $i=1,2,\dots,n$, let $\Gamma_i$ be a signed graph
with exactly two distinct eigenvalues $\pm \theta_i$, where $\Gamma_1,\dots,\Gamma_{n-1}$ are bipartite.
Then $\Gamma_{\widetilde{\Box},L}^n$, $\Gamma_{\widetilde{\Box},R}^n$, $\Gamma_{\widetilde{\bowtie},L}^n$ and $\Gamma_{\widetilde{\bowtie},R}^n$
have exactly two distinct eigenvalues $\pm \sqrt{\sum_{i=1}^n\theta_i^2}$, $\pm \sqrt{\sum_{i=1}^n\theta_i^2}$,
$\pm \sqrt{\theta_n^2\prod_{i=1}^{n-1}(\theta_i^2+1)}$
and $\pm \sqrt{\sum_{k=1}^n\prod_{i=k}^{n}\theta_i^2}$, respectively.
\end{cor}

\begin{exa}\label{example-Qn}
{\rm
Let $\Gamma_i=(K_2, +)$ for each $i=1,2,\dots,n$. Then

(i) each of $\Gamma_{\widetilde{\Box},R}^n$, $\Gamma_{\widetilde{\Box},L}^n$ and $\Gamma_{\widetilde{\bowtie},R}^n$
is a signed graph of $Q_n$ whose eigenvalues are $\pm \sqrt{n}$;

(ii) $\Gamma_{\widetilde{\bowtie},L}^n$ is a signed graph of $K_{2^{n-1},2^{n-1}}$
and its eigenvalues are $\pm \sqrt{2^{n-1}}$.
}
\end{exa}

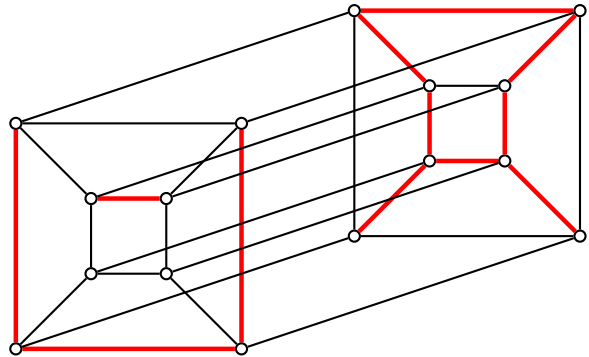
\begin{figure}[ht]
\begin{center}
\begin{pspicture}(-3,-6.5)(3,-2)

\cnode(-3,-3.5){2.5pt}{A1} \cnode(0,-3.5){2.5pt}{A2} \cnode(0,-6.5){2.5pt}{A3} \cnode(-3,-6.5){2.5pt}{A4}
\ncline{A1}{A2}\ncline[linewidth=0.6mm,linecolor=red]{A1}{A4}\ncline[linewidth=0.6mm,linecolor=red]{A3}{A2}\ncline[linewidth=0.6mm,linecolor=red]{A3}{A4}
\cnode(-2,-4.5){2.5pt}{a1} \cnode(-1,-4.5){2.5pt}{a2} \cnode(-1,-5.5){2.5pt}{a3} \cnode(-2,-5.5){2.5pt}{a4}
\ncline[linewidth=0.6mm,linecolor=red]{a1}{a2}\ncline{a1}{a4}\ncline{a3}{a2}\ncline{a3}{a4}
\ncline{a1}{A1}\ncline{a2}{A2}\ncline{a3}{A3}\ncline{a4}{A4}

\cnode(1.5,-2){2.5pt}{B1} \cnode(4.5,-2){2.5pt}{B2} \cnode(4.5,-5){2.5pt}{B3} \cnode(1.5,-5){2.5pt}{B4}
\ncline[linewidth=0.6mm,linecolor=red]{B1}{B2}\ncline{B1}{B4}\ncline{B3}{B2}\ncline{B3}{B4}
\cnode(2.5,-3){2.5pt}{b1} \cnode(3.5,-3){2.5pt}{b2} \cnode(3.5,-4){2.5pt}{b3} \cnode(2.5,-4){2.5pt}{b4}
\ncline{b1}{b2}\ncline[linewidth=0.6mm,linecolor=red]{b1}{b4}\ncline[linewidth=0.6mm,linecolor=red]{b3}{b2}\ncline[linewidth=0.6mm,linecolor=red]{b3}{b4}
\ncline[linewidth=0.6mm,linecolor=red]{b1}{B1}\ncline[linewidth=0.6mm,linecolor=red]{b2}{B2}
\ncline[linewidth=0.6mm,linecolor=red]{b3}{B3}\ncline[linewidth=0.6mm,linecolor=red]{b4}{B4}

\ncline{a1}{b1}\ncline{a2}{b2}\ncline{a3}{b3}\ncline{a4}{b4}
\ncline{A1}{B1}\ncline{A2}{B2}\ncline{A3}{B3}\ncline{A4}{B4}

\end{pspicture}
\caption{{\small The signed graph $\Gamma_{\widetilde{\Box},L}^4$ of $Q_4$ in Example \ref{example-Qn}.}\label{fig2}}
\end{center}
\end{figure}

\begin{exa}\label{example-Q2n}
{\rm
For each $i=1,2,\dots,n$, let $\Gamma_i=(K_{2,2},\sigma)$ be
the signed graph of $K_{2,2}$ with exactly one negative edge. Then

(i) $\Gamma_{\widetilde{\Box},R}^n$ and $\Gamma_{\widetilde{\Box},L}^n$
are signed graphs of $Q_{2n}$ whose eigenvalues are
$\pm \sqrt{2n}$; 

(ii) the eigenvalues of $\Gamma_{\widetilde{\bowtie},L}^n$
are $\pm \sqrt{2\cdot3^{n-1}}$; 

(iii) the eigenvalues of $\Gamma_{\widetilde{\bowtie},R}^n$
are $\pm \sqrt{\sum_{k=1}^n2^k}= \pm \sqrt{2^{n+1}-2}$. 
}
\end{exa}

By Definition \ref{Def-product-n-graphs}, Theorem \ref{thm-symmetric} can be generalized to Theorem \ref{thm-sym-n}.

\begin{thm}\label{thm-sym-n}
For $n\geq 2$ and $i=1,2,\dots,n-1$, let $\Gamma_i$ be a signed bipartite graph
and $\Gamma_n$ be a signed graph.

(i) The spectrum of $\Gamma_{\widetilde{\bowtie},R}^n$ is symmetric if and only
if $\Gamma_{n-1}$ is balanced or the spectrum of $\Gamma_{n}$ is symmetric.

(ii) The spectrum of every graph in $\{\Gamma_{\widetilde{\Box},R}^n, \Gamma_{\widetilde{\Box},L}^n,
\Gamma_{\widetilde{\bowtie},L}^n \}$ is symmetric if and only if
there exists an integer $i\in\{1,\dots,n-1\}$ such that $\Gamma_i$ is balanced
or the spectrum of $\Gamma_{n}$ is symmetric.
\end{thm}

\begin{pf}
We only need to consider $n\geq 3$.

(i) By Theorem \ref{thm-symmetric}, the spectrum of $\Gamma_{\widetilde{\bowtie},R}^n
=\Gamma_{\widetilde{\bowtie},R}^{n-1}\widetilde{\bowtie}\Gamma_{n}$
is symmetric if and only if $\Gamma_{\widetilde{\bowtie},R}^{n-1}$ is balanced or
the spectrum of $\Gamma_{n}$ is symmetric.
Since the semi-strong product of a graph $G$ and a bipartite graph $H$ is
balanced if and only if $H$ is balanced, we have
$\Gamma_{\widetilde{\bowtie},R}^{n-1}
=\Gamma_{\widetilde{\bowtie},R}^{n-2}\widetilde{\bowtie} \Gamma_{n-1}$
is balanced if and only if $\Gamma_{n-1}$ is balanced. Thus, (i) is proved.

(ii) By Theorem \ref{thm-symmetric}, the spectrum of $\Gamma_{\widetilde{\Box},R}^n=
\Gamma_{\widetilde{\Box},R}^{n-1}\widetilde{\Box} \Gamma_n$ is symmetric if and only if
$\Gamma_{\widetilde{\Box},R}^{n-1}$ is balanced or the spectrum of $\Gamma_{n}$ is symmetric.
Since $\Gamma_{\widetilde{\Box},R}^{n-1}$ is balanced if and only if
there exists an integer $i\in\{1,\dots,n-1\}$ such that $\Gamma_i$ is balanced.
So the conclusion for $\Gamma_{\widetilde{\Box},R}^n$ is proved.

By Theorem \ref{thm-symmetric}, the conclusion holds for
$\Gamma_{\widetilde{\Box},L}^2=\Gamma_{n-1}\widetilde{\Box}\Gamma_n$
(resp. $\Gamma_{\widetilde{\bowtie},L}^2=\Gamma_{n-1}\widetilde{\bowtie}\Gamma_n$). By induction on $n$,
assume that the spectrum of $\Gamma_{\widetilde{\Box},L}^{n-1}$
(resp. $\Gamma_{\widetilde{\bowtie},L}^{n-1}$) is symmetric if and only if
there exists an integer $i\in\{2,\dots,n-1\}$ such that $\Gamma_i$ is balanced
or the spectrum of $\Gamma_{n}$ is symmetric.
By Theorem \ref{thm-symmetric}, the spectrum of $\Gamma_{\widetilde{\Box},L}^n=
\Gamma_1\widetilde{\Box}\Gamma_{\widetilde{\Box},L}^{n-1}$
(resp. $\Gamma_{\widetilde{\bowtie},L}^n=
\Gamma_1\widetilde{\bowtie}\Gamma_{\widetilde{\bowtie},L}^{n-1}$) is symmetric if and only if
$\Gamma_1$ is balanced or the spectrum of $\Gamma_{\widetilde{\Box},L}^{n-1}$
(resp. $\Gamma_{\widetilde{\bowtie},L}^{n-1}$) is symmetric.
By induction,  the conclusion for $\Gamma_{\widetilde{\Box},L}^n$
(resp. $\Gamma_{\widetilde{\bowtie},L}^n$) is proved.
\end{pf}

\bigskip
Now, Theorem \ref{thm-Delta-2-min} is generalized to the following theorem.

\begin{thm}\label{thm-Delta-n-graph}
For $i=1,2,\dots,n$, let $\Gamma_i=(G_i,\sigma_i)$ be a signed graph of order $N_i$
and $\theta_i^2$ be the minimum eigenvalue of $A(\Gamma_i)^2$, where $G_1,\dots,G_{n-1}$ are bipartite. 
Let $H_{\Box}$, $H_{\bowtie,L}$ and $H_{\bowtie,R}$ be any $(\lfloor\frac{1}{2}\prod_{i=1}^nN_i\rfloor+1)$-vertex induced
subgraph of $\Gamma_{\widetilde{\Box},L}^n$, $\Gamma_{\widetilde{\bowtie},L}^n$ and $\Gamma_{\widetilde{\bowtie},R}^n$, respectively.

(i) If there exists an integer $i\in\{1,2,\dots,n-1\}$ such that $\Gamma_i$ is balanced
or the spectrum of $\Gamma_n$ is symmetric, then
$\Delta(H_{\Box})\geq \sqrt{\sum_{i=1}^n\theta_i^2}$ and
$\Delta(H_{\bowtie,L})\geq \sqrt{\theta_n^2\prod_{i=1}^{n-1}(\theta_i^2+1)}$;

(ii) If $\Gamma_{n-1}$ is balanced or the spectrum of $\Gamma_n$ is symmetric, then
$\Delta(H_{\bowtie,R})\geq \sqrt{\sum_{k=1}^n\prod_{i=k}^{n}\theta_i^2}$.
\end{thm}

\begin{pf}
For simplicity, let $N=\prod_{i=1}^nN_i$. 
Since $H_{\Box}$, $H_{\bowtie,L}$ and $H_{\bowtie,R}$ are $(\lfloor\frac{N}{2}\rfloor+1)$-vertex induced subgraphs of
$\Gamma_{\widetilde{\Box},L}^n$, $\Gamma_{\widetilde{\bowtie},L}^n$ and
$\Gamma_{\widetilde{\bowtie},R}^n$, respectively.
By Lemma \ref{lem-Delta-lambda},
$$
\Delta(H_{\Box})\geq \lambda_{\lceil\frac{1}{2}N\rceil}(\Gamma_{\widetilde{\Box},L}^n),
 \ \Delta(H_{\bowtie,L})\geq \lambda_{\lceil\frac{1}{2}N\rceil}(\Gamma_{\widetilde{\bowtie},L}^n)\  \text{and} \
\Delta(H_{\bowtie,R})\geq \lambda_{\lceil\frac{1}{2}N\rceil}(\Gamma_{\widetilde{\bowtie},R}^n).
$$
By Theorem \ref{thm-eigenvalues-n-graph}, the minimum eigenvalues of
$A(\Gamma_{\widetilde{\Box},L}^n)^2$, $A(\Gamma_{\widetilde{\bowtie},L}^n)^2$
and $A(\Gamma_{\widetilde{\bowtie},R}^n)^2$ are obtained. Thus, by
Theorem \ref{thm-sym-n}, the spectrums of $\Gamma_{\widetilde{\Box},L}^n$,
$\Gamma_{\widetilde{\bowtie},L}^n$ and $\Gamma_{\widetilde{\bowtie},R}^n$ are symmetric and so
$\lambda_{\lceil\frac{1}{2}N\rceil}(\Gamma_{\widetilde{\Box},L}^n)$ $= \sqrt{\sum\nolimits_{i=1}^n\theta_i^2}$,
$\lambda_{\lceil\frac{1}{2}N\rceil}(\Gamma_{\widetilde{\bowtie},L}^n)=\sqrt{\theta_n^2\prod\nolimits_{i=1}^{n-1}(\theta_i^2+1)}$ and
$\lambda_{\lceil\frac{1}{2}N\rceil}(\Gamma_{\widetilde{\bowtie},R}^n)=\sqrt{\sum\nolimits_{k=1}^n\prod\nolimits_{i=k}^{n}\theta_i^2}.$
Combining these (in)equalities, the results follow.
\end{pf}

\begin{cor}\label{cor-DeltaH-n}
For $i=1,2,\dots,n$, let $\Gamma_i=(G_i,\sigma_i)$ be a signed graph of order $N_i$
with exactly two distinct eigenvalues $\pm\theta_i$, where $G_1,\dots,G_{n-1}$ are bipartite. 
Let $H_{\Box}$, $H_{\bowtie,L}$ and $H_{\bowtie,R}$ be any $(\lfloor\frac{1}{2}\prod_{i=1}^nN_i\rfloor+1)$-vertex induced
subgraph of $\Gamma_{\widetilde{\Box},L}^n$, $\Gamma_{\widetilde{\bowtie},L}^n$ and $\Gamma_{\widetilde{\bowtie},R}^n$, respectively.
Then
$\Delta(H_{\Box})\geq \sqrt{\sum_{i=1}^n\theta_i^2}$,
$\Delta(H_{\bowtie,L})\geq \sqrt{\theta_n^2\prod_{i=1}^{n-1}(\theta_i^2+1)}$ and
$\Delta(H_{\bowtie,R})\geq \sqrt{\sum_{k=1}^n\prod_{i=k}^{n}\theta_i^2}$.
\end{cor}

When $\Gamma_i=(K_2,+)$ for each $i=1,2,\dots,n$ in Corollary \ref{cor-DeltaH-n},
$\Gamma_{\widetilde{\Box},L}^n$ is the signed graph of hypercube $Q_n$.
Therefore, Corollary \ref{cor-DeltaH-n} implies Huang's theorem.

\begin{exa}{\rm
For $i=1,2,\dots,n$, let $\Gamma_i=(K_{2^t,2^t},\sigma)$ be the signed graph $K_{2^t,2^t}$
with exactly two distinct eigenvalues $\pm \sqrt{2^t}$.
For any integer $n\geq 1$ and $t\geq 0$, let $H_{\Box}$, $H_{\bowtie,L}$ and $H_{\bowtie,R}$
be any $(2^{n(t+1)-1}+1)$-vertex induced subgraph of $\Gamma_{\widetilde{\Box},L}^n$,
$\Gamma_{\widetilde{\bowtie},L}^n$ and $\Gamma_{\widetilde{\bowtie},R}^n$ respectively.
Then
$\Delta(H_{\Box})\geq \sqrt{2^t \cdot n}$,
$\Delta(H_{\bowtie,L})\geq \sqrt{2^t(2^t+1)^{n-1}}$ and
$\Delta(H_{\bowtie,R})\geq \sqrt{\sum_{k=1}^n2^{kt}}$.
}
\end{exa}

\section{Concluding remarks}

{\bf I.} Corollary \ref{cor-kron-eig}, Corollary \ref{cor-eig-G1G2} and Corollary \ref{cor-eig-G1Gn}
provide product methods to construct signed graphs with exactly
two distinct eigenvalues of opposite signatures from factor graph $\Gamma_1$ and $\Gamma_2$.
There are many options for the factor graph, such as
the signed graphs of $Q_n$ and $K_{2^n,2^n}$ in Example \ref{example-Qn},
$T_{2n}$ in Lemma \ref{lem-T2n},
$S_{14}$ in Lemma \ref{lem-S14},
the signed graph of $K_n$ in Example \ref{exa-Kn},
signed graphs in Examples \ref{example-complete-k-part-graph},
\ref{example-Hn-AGamma}, \ref{example-Q2n} and so on.

\bigskip
{\bf II.}  If the following conjecture is true, it would provide a way to construct an infinite
family of $d$-regular Ramanujan graphs by $2$-lift of graphs.

\begin{con}(Bilu-Linial \cite{BiLi06})\label{con-BiLi06}
Every connected $d$-regular graph $G$ has a signature $\sigma$ such that 
$\rho(G,\sigma)\leq 2\sqrt{d-1}$.
\end{con}

Gregory considered the following Conjecture \ref{con-Greg12} without the regularity assumption on $G$.

\begin{con}(Gregory \cite{Greg12})\label{con-Greg12}
If $G$ is a nontrivial graph with maximum degree $\Delta>1$, then there exists
a signed graph $\Gamma=(G,\sigma)$ such that $\rho(\Gamma)\leq 2\sqrt{\Delta-1}$.
\end{con}

By Theorem \ref{thm-eigenvalues} (i), we have the following theorem.

\begin{thm}\label{thm-Ramanujan}
For $i=1,2$, let $G_i$ be a graph with maximum degree $\Delta_i$ and
$\Gamma_i=(G_i,\sigma_i)$ be a signed graph such that
$\rho(\Gamma_i)\leq 2\sqrt{\Delta_i-1}$. If $G_1$ is bipartite, then
$$\rho(\Gamma_1\widetilde{\Box} \Gamma_2)\leq 2\sqrt{\Delta_1+\Delta_2-2}.$$
\end{thm}

Since $\rho(\Gamma_1\widetilde{\Box} \Gamma_2)=\sqrt{\rho(\Gamma_1)^2+\rho(\Gamma_2)^2}$
and $\Delta(\Gamma_1\Box \Gamma_2)=\Delta(\Gamma_1)+\Delta(\Gamma_2)$,
Theorem \ref{thm-Ramanujan} shows that if Conjecture \ref{con-Greg12} holds for $\Gamma_1$ and $\Gamma_2$,
then Conjecture \ref{con-Greg12} also holds for the signed Cartesian product of them.

\bigskip
{\bf III.} The method which is utilized to construct a larger weighing matrix can
construct a larger signed graph with exactly two distinct eigenvalues $\pm\theta$ from small graphs.
Conversely, the ideas of signed Cartesian product and semi-strong product in our paper
can also be applied to construct a weighing matrix.
If for $i=1,2$, $W_i$ is a weighing matrix of order $n_i$ and weight $k_i$, then we can construct
weighing matrices as follows
\begin{equation*}
W(4n_1n_2, k_1+k_2) = \left[
\begin{array}{cc}
O_{n_1} & W_1 \\
W_1^T & O_{n_1}
\end{array}
\right] \otimes I_{2n_2} +
\left[
\begin{array}{cc}
I_{n_1} & O_{n_1} \\
O_{n_1} & -I_{n_1}
\end{array}
\right] \otimes \left[
\begin{array}{cc}
O_{n_2} & W_2 \\
W_2^T & O_{n_2}
\end{array}
\right],
\end{equation*}
\begin{equation*}
W(4n_1n_2, (k_1+1)k_2) =  \left[
\begin{array}{cc}
I_{n_1} & W_1  \\
W_1^T  & -I_{n_1}
\end{array}
\right] \otimes \left[
\begin{array}{cc}
O_{n_2} & W_2 \\
W_2^T & O_{n_2}
\end{array}
\right],
\end{equation*}
\begin{equation*}
W(2n_1n_2, (k_1+1)k_2) = \left[
\begin{array}{cc}
I_{n_1} & W_1  \\
W_1^T  & -I_{n_1}
\end{array}
\right] \otimes W_2 .
\end{equation*}
Furthermore, if $W_2$ is symmetric, then we can construct weighing matrix
\begin{equation*}
W(2n_1n_2, k_1+k_2) = \left[
\begin{array}{cc}
I_{n_1}\otimes W_2 & W_1\otimes I_{n_2} \\
W_1^T\otimes I_{n_2} & -I_{n_1}\otimes W_2
\end{array}
\right].
\end{equation*}

\section*{Acknowledgement}

The research of Zhen-Mu Hong is supported by NNSFC (No. 11601002),
Outstanding Young Talents International Visiting Program of Anhui Provincial Department of Education (No. gxgwfx2018031)
and Key Projects in Natural Science Research of Anhui Provincial Department of Education (No. KJ2016A003).
The research of Hong-Jian Lai is supported by NNSFC (Nos. 11771039 and 11771443).

\end{document}